\documentclass{article}
   \usepackage{amssymb} 
 \usepackage[totalwidth=410pt, totalheight=620pt]{geometry}
\newcommand{\g}{\mbox{$\bf g$}}

\newcommand{\h}{\mbox{\textbf{h}}}
\newcommand{\n}{\mbox{\textbf{n}}}

\newcommand{\al}{\alpha}

\newcommand{\la}{\lambda}
\newcommand{\La}{\Lambda}

\newcommand{\Th}{\Theta}
\newcommand{\gt}{\theta}
\newcommand{\mb}{\mbox}

\newcommand{\Mklz}[2]{\left\{\left.\;#1\;\right|\; #2\;\right\}}
\newcommand{\nW}{\mbox{$\Delta^-$}}
\newcommand{\pW}{\mbox{$\Delta^+$}}
\newcommand{\rW}{\mbox{$\Delta_{re}$}}
\newcommand{\nrW}{\mbox{$\Delta_{re}^-$}}
\newcommand{\prW}{\mbox{$\Delta_{re}^+$}}

\newcommand{\F}{\mathbb{F}}

\newcommand{\N}{\mathbb{N}}
\newcommand{\Nn}{\mathbb{N}_0}
\newcommand{\Q}{\mathbb{Q}}

\newcommand{\R}{\mathbb{R}}

\newcommand{\Z}{\mathbb{Z}}
\newcommand{\We}{\mbox{$\mathcal W$}}

\newcommand{\iB}[2]{\left(#1\mid#2\right)}

\newcommand{\TD}{\mbox{$\widehat{T}$}}

\newcommand{\GD}{\mbox{$\widehat{G}$}}
\newcommand{\GfD}{\mbox{$\widehat{G_f}$}}

\newcommand{\FK}[1]{\mbox{$\F\,[#1]$}}

\newcommand{\Fi}[1]{{\mathcal F}_{#1}}

\newcommand{\Proof}{\mbox{\bf Proof: }}
\newcommand{\qed}{\mb{}\hfill\mb{$\square$}\\}

\newcommand{\Spm}{\mbox{Specm\,}}

\newcommand{\ti}{\tilde}
\newcommand{\res}[1]{\!\mid_{#1}}

\newcommand{\ro}{{\,\bf \triangleright\,}}
\newcommand{\lo}{{\,\bf \triangleleft\,}}
\begin{document}
\newtheorem{Theorem}{Theorem}[section]
\newtheorem{Definition}[Theorem]{Definition}
\newtheorem{Proposition}[Theorem]{Proposition}
\newtheorem{Corollary}[Theorem]{Corollary}
\newtheorem{Conjecture}[Theorem]{Conjecture}
\newtheorem{Remark}[Theorem]{Remark}
\newtheorem{Remarks}[Theorem]{Remarks}
\newtheorem{Lemma}[Theorem]{Lemma}
\newtheorem{Example}[Theorem]{Example}
\newtheorem{Examples}[Theorem]{Examples}
%
%%%%%%%%%%%%%%%%%%%%%%%%%%%%%%%%%%%%%%%%%%%%%%%%%%%%%%%%%%%%%%%%%%%%%%%%%%%%%%%%%%%%%%%%%%%%%%%%%%%%%%%%%%%%%%%%%%%%%%%%%%%%%%%%%%%%%%%%%%%%%
%
\title{Integrating infinite-dimensional Lie algebras by a Tannaka
  reconstruction (Part II)}
\author{Claus Mokler\thanks{Supported by the Deutsche Forschungsgemeinschaft.}\\\\ Universit\"at Wuppertal, Fachbereich C - Mathematik\\  Gau\ss stra\ss e 20\\ D-42097 Wuppertal, Germany\vspace*{1ex}\\ 
          mokler@math.uni-wuppertal.de}
\date{}
\maketitle
\begin{abstract}\noindent 

Let $\g$ be a Lie algebra over a field $\F$ of characteristic zero, let $\mathcal C$ be a certain tensor category of representations of $\g$, 
and ${\mathcal C}^{du}$ a certain category of duals. In \cite{M4} we associated to $\mathcal C$ and ${\mathcal C}^{du}$ 
by a Tannaka reconstruction a monoid $M$ with a coordinate ring $\FK{M}$ of matrix coefficients, as well as a Lie algebra $Lie(M)$. 
We interpreted the Tannaka monoid $M$ algebraic geometrically as a weak algebraic monoid with Lie algebra $Lie(M)$. The monoid $M$ acts 
by morphisms of varieties on every object $V$ of $\mathcal C$. The Lie algebra $Lie(M)$ acts on $V$ by the differentiated action. 

In particular, we showed in \cite{M4}: If the Lie algebra $\g$ is generated by by one-parameter elements, then it identifies in a natural way 
with a subalgebra of $Lie(M)$, and there exists a subgroup of the unit group of $M$, which is dense in $M$. 
In the present paper we introduce the coordinate ring of regular functions on this dense subgroup of $M$, as well as the algebra of linear regular 
functions on the universal enveloping algebra $U(\g)$ of $\g$, and investigate their relation. We investigate and describe various coordinate 
rings of matrix coefficients associated to categories of integrable representations of $\g$. We specialize to
integrable representations of Kac-Moody algebras and free Lie algebras. 
Some results on coordinate rings of Kac-Moody groups obtained by V. G. Kac and D. Peterson, some coordinate rings of the associated
groups of linear algebraic integrable Lie algebras defined by V. G. Kac, and some results
on coordinate rings of free Kac-Moody groups obtained by Y. Billig and A. Pianzola fit into this context.

We determine the Tannaka monoid associated to the full subcategory of integrable
representations in the category $\mathcal O$ of a Kac-Moody algebra and to its category of full
duals. Its Zariski-open dense unit group is the formal Kac-Moody group. We give various descriptions of its coordinate ring of matrix
coefficients. We show that its Lie algebra is the formal Kac-Moody algebra.
\end{abstract}
\vspace*{1ex} 
\section*{Introduction}
Finite-dimensional Lie algebras arise as an important aid in the investigation of fi\-nite-\-di\-men\-sion\-al Lie groups and algebraic groups, 
encoding their local structure. 
In the infinite-dimensional situation one often finds natural infinite-dimensional Lie algebras, which also have a rich representation
theory. But the Lie algebras are defined and obtained without the help of any groups. One faces the problem to ``integrate infinite-dimensional 
Lie algebras'', i.e., to construct groups associated in an appropriate way to these Lie algebras. Important examples are Kac-Moody algebras 
and their associated Kac-Moody groups.

In the paper \cite{M4} we gave an approach to the integration of infinite-dimensional Lie algebras based on the Tannaka reconstruction.  
Starting with a Lie algebra $\g$ over a field $\F$ of characteristic zero, a certain category $\mathcal C$ of representations 
of $\g$, and a certain category of duals ${\mathcal C}^{du}$, we associated to $\mathcal C$ and ${\mathcal C}^{du}$ the Tannaka monoid $M$ with 
coordinate ring of matrix coefficients $\FK{M}$, as well as a Lie algebra $Lie(M)$. We interpreted the Tannaka monoid $M$ algebraic geometrically as an 
irreducible  weak algebraic monoid with Lie algebra $Lie(M)$. The monoid $M$ acts 
by morphisms of varieties on every object $V$ of $\mathcal C$. The Lie algebra $Lie(M)$ acts on $V$ by the differentiated action. 
The pair of categories $\mathcal C$ and ${\mathcal C}^{du}$ is good for integrating $\g$, if the Lie algebra $\g$ identifies in a natural 
way with a Lie subalgebra of $Lie(M)$. We investigated the Tannaka monoid $M$, its coordinate ring of matrix coefficients $\FK{M}$, 
and the Lie algebra $Lie(M)$ in this situation.
In general, it is a difficult problem to decide if some categories $\mathcal C$ and ${\mathcal C}^{du}$ are good for integrating $\g$. We showed 
that if the Lie algebra $\g$ is generated by integrable locally finite elements, then $\mathcal C$ and $\mathcal C^{du}$ are good for integrating
$\g$.

One-parameter elements are particular examples of integrable locally finite elements. If $\g$ is generated by one-parameter elements, 
then there also exists a subgroup of the unit group of $M$, which is dense in $M$. In Section \ref{regfct} of the present paper we introduce 
the coordinate ring of regular functions on this dense subgroup of $M$.
For particular examples of Lie algebras, categories of representations, and generating sets of one-parameter elements, this subgroup and its 
coordinate ring specialize to the (minimal) Kac-Moody groups and its coordinate rings of regular functions defined and investigated 
by V. G. Kac and D. Peterson in \cite{KP1}, \cite{KP2}, 
to the associated groups of linear algebraic integrable Lie algebras and its coordinate rings of regular functions defined by V. G. Kac in \cite{K1}, and
to the free Kac-Moody groups and its coordinate rings of polynomial functions defined and investigated by Y. Billig and A. Pianzola in \cite{BiPi}.
We introduce the counterpart of the coordinate ring of regular functions, the algebra of linear regular functions on the universal enveloping algebra 
$U(\g)$ of $\g$, with which it is easier to work algebraically. We give a theorem, which describes when both algebras are isomorphic. 
In particular, applied to the algebra of polynomial functions on the free Kac-Moody group, giving this isomorphism is roughly equivalent 
to a description of the algebra of polynomial functions on the free Kac-Moody group obtained Y. Billig and A. Pianzola in \cite{BiPi}.

In Section \ref{mop} we define various coordinate rings of matrix coefficients associated to categories of integrable representations 
of a Lie algebra $\g$. This is motivated by the definition of some coordinate rings of the associated groups of linear algebraic integrable Lie algebras 
in \cite{K1}, and by some results in relation to the free Kac-Moody group and the shuffle algebra in \cite{BiPi}. 
Depending on certain assumptions which have to be satisfied we characterize the corresponding subalgebras of the coordinate rings of
regular functions, as well as the corresponding subalgebras of the dual of the
universal enveloping algebra $U(\g)$ of $\g$.
We specialize to integrable representations of free Lie algebras and symmetrizable Kac-Moody algebras. 
We also show how some of the results related to the shuffle algebra obtained in \cite{BiPi} and \cite{Pi2} can be interpreted over a field of
characteristic zero by the Tannaka reconstruction.

Even if we know that $\mathcal C$ and ${\mathcal C}^{du}$ are good for integrating an infinite-dimensional Lie algebra, it may be quite hard, and a
long way to determine the associated Tannaka monoid and its Lie algebra. An
example has been treated in \cite{M1}. The full subcategory ${\mathcal O}_{int}$ of the category
$\mathcal O$ of a symmetrizable Kac-Moody algebra $\g$, whose objects are the integrable 
$\g$-modules contained in $\mathcal O$, is one possible generalization of the category of 
finite-dimensional representations of a semisimple Lie algebra. It keeps the complete reducibility theorem, i.e., every integrable 
$\g$-module contained in $\mathcal O$ is a sum of integrable irreducible highest weight modules. In \cite{M1} we determined the Tannaka monoid 
associated to this category and its category of integrable duals. Its Zariski-open dense unit group is the (minimal) Kac-Moody group. We showed that its 
Lie algebra identifies with the Kac-Moody algebra. 
In Section \ref{KMintfull} we determine the Tannaka monoid associated to the category of integrable $\g$-modules contained in $\mathcal O$ and its 
category of full duals. Its Zariski-open dense unit group is the formal Kac-Moody group. We give several descriptions of its coordinate ring of matrix 
coefficients. We show that its Lie algebra identifies with the formal Kac-Moody algebra. It is possible to do this in a short way by making use of some 
results of \cite{M2}.

Many examples remain to be treated. In particular, for a symmetrizable Kac-Moody algebra $\g$ it would be important to determine the Tannaka monoid 
associated to the category of integrable $\g$-modules with point separating integrable duals and to its category of
integrable duals. We conjecture that the Tannaka monoid is the (minimal)
Kac-Moody group, and that its Lie algebra identifies with the Kac-Moody algebra.
%
%%%%%%%%%%%%%%%%%%%%%%%%%%%%%%%%%%%%%%%%%%%%%%%%%%%%%%%%%%%%%%%%%%%%%%%%%%%%%%%%%%%%%%%%%%%%%%%%%%%%%%%%%%%%%%%%%%%%%%%%%%%%%%%%%%%%%%%%%%%%%%%%%%%%%
\tableofcontents
%%%%%%%%%%%%%%%%%%%%%%%%%%%%%%%%%%%%%%%%%%%%%%%%%%%%%%%%%%%%%%%%%%%%%%%%%%%%%%%%%%%%%%%%%%%%%%%%%%%%%%%%%%%%%%%%%%%%%%%%%%%%%%%%%%%%%%%%%%%%%%%%%%%%%%
%
%\newpage
\section{Preliminaries}
This paper continues the paper \cite{M4}. We use the notation introduced in \cite{M4}. We also refer
for the basic definitions to \cite{M4}. Before starting to read Section \ref{regfct} and Section \ref{mop} of this paper, please review in particular Section 2.1, and
Section 2.4 of \cite{M4}. Before starting to read Section \ref{KMintfull} of this paper please review Section 3 of \cite{M4}.
%
%%%%%%%%%%%%%%%%%%%%%%%%%%%%%%%%%%%%%%%%%%%%%%%%%%%%%%%%%%%%%%%%%%%%%%%%%%%%%%%%%%%%%%%%%%%%%%%%%%%%%%%%%%%%%%%%%%%%%%%%%%%%%%%%%%%%%%%%%%%%%%%%%% 
%
%
\section{Regular functions\label{regfct}}
%
%%%%%%%%%%%%%%%%%%%%%%%%%%%%%%%%%%%%%%%%%%%%%%%%%%%%%%%%%%%%%%%%%%%%%%%%%%%%%%%%%%%%%%%%%%%%%%%%%%%%%%%%%%%%%%%%%%%%%%%%%%%%%%%%%%%%%%%%%%%%%%%%%%
%
%
%
Let $\g$ be a Lie algebra over a field $\F$ of characteristic 0. Fix a category $\mathcal C$, and a category ${\mathcal C}^{du}$ of duals as in 
Section 2.1 of \cite{M4}. Let $M$ be the associated Tannaka monoid equipped with its coordinate ring $\FK{M}$ of matrix coefficients as defined 
and described in Section 2.1 of \cite{M4}. Let $Lie(M)$ be its Lie algebra, which we realize as in Section 2.1 of \cite{M4} as a Lie subalgebra 
of the Lie algebra $Nat$. Also as described in Section 2.1 of \cite{M4} we identify $\g$ with a subalgebra of $Nat$.

Regular functions for Kac-Moody groups have been introduced in \cite{KP2}. Regular functions for groups associated to linear algebraic
integrable Lie algebras in \cite{K1}, \S 1.8. Polynomial functions for groups called free Kac-Moody groups have been introduced and 
investigated in \cite{BiPi}.
The concept of regular functions fits into the context of the Tannaka reconstruction of \cite{M4}, including these cases. 
For a simple presentation we restrict here to the case where the Lie algebra $\g$ is generated by one-parameter elements. More generally, 
it is possible to define regular functions if the Lie algebra $\g$ is generated by integrable locally finite elements.\vspace*{1ex}  

Recall the following parts of Theorem 2.35, Remarks 2.36, and Theorem 2.37, Remarks 2.38 of \cite{M4}, which are important
for our considerations, and which also introduce some of the notation we use in this paper.

\begin{Theorem} \label{one1II} Let $e\in \g\setminus\{0\}$ such that for every $\g$-module $V$ contained in $\mathcal C$ 
\begin{itemize}
\item[(1)] $e_V$ is locally nilpotent,
\item[(2)] $\exp(t e_V)\in End_{V^{du}}(V)$ for all $t\in\F$.
\end{itemize}
Then for every $t\in \F$ there exists an element $\exp(te)\in M$, which acts on every $\g$-module $V$ contained 
in $\mathcal C$ by $\exp(t e_V)$.
Equip $\F$ with the coordinate ring of polynomial functions $\FK{t}$. The map
\begin{eqnarray*}
          \kappa_e:\;   (\F,+)   &\to      & \;\;\;M\\
          t \;\;\;\;  & \mapsto & \exp(te) 
\end{eqnarray*}
is a morphism of monoids, which is also an embedding of sets with coordinate rings. Its image $V_e$ is closed in $M$, and $Lie(V_e)=\F e$.
\end{Theorem}
\begin{Remark}\label{Rone1II} Condition (2) of the theorem is satisfied if $V^{du}=V^*$, or if $(e_V)^{du}$ is a locally nilpotent endomorphism of $V^{du}$.  
\end{Remark}

\begin{Theorem}\label{one2II}
Let $e\in\g$ such that for every $\g$-module $V$ contained in $\mathcal C$ 
\begin{itemize}
\item[(1)] $e_V$ is diagonalizable with integer eigenvalues,
\item[(2)] $t^{e_V}\in End_{V^{du}}(V)$ for every $t\in\F^\times$. 
\end{itemize}
Then for every $t\in\F^\times$ there exists an element $t^e\in M$, which acts on every module $V$ contained in 
$\mathcal C$ by $t^{e_V}$. Equip $\F^\times$ with the coordinate ring of Laurent polynomial
  functions $\FK{t, t^{-1}}$. Let $\Z_e$ be the set of eigenvalues of $e_V$ for all objects $V$ of $\mathcal C$. Then
\begin{eqnarray*}
  \kappa_e:\;(\F^\times ,\cdot)  &\to      & M\\
   t\;\;\;\;  & \mapsto &  t^e
\end{eqnarray*}
is a morphism of monoids, which is also a morphism of sets with coordinate rings. 
$\Z_e$ is a submonoid of $(\Z,+)$, and the image of the comorphism $\kappa_e^*$ is given by
\begin{eqnarray*}
   \kappa_e^*(\FK{M})=\bigoplus_{b\in \Z_e} \F \,t^b\;.
\end{eqnarray*}
Denote the image of the map $\kappa_e$ by $V_e$. We have $Lie(V_e)=\F e$.  
\end{Theorem}
\begin{Remarks}\label{Rone2II} (1) Condition (2) of the theorem is satisfied if $V^{du}=V^*$, or if $(e_V)^{du}$ is a diagonalizable 
endomorphism of $V^{du}$.

(2) If $e$ acts by positive and negative eigenvalues on the objects of $\mathcal C$ then $\Z_e=a\Z$, where $a$ is the smallest absolute value of the 
non-zero eigenvalues of $e_V$ for all objects $V$ of $\mathcal C$.. 
\end{Remarks}

We call elements $e\in\g\setminus\{0\}$ as used in Theorem \ref{one1II} {\it locally nilpotent one-parameter elements}.
We call elements $e\in\g\setminus\{0\}$ as used in Theorem \ref{one2II} {\it diagonalizable one-parameter elements}. 
A diagonalizable one-parameter element $e$ is called $(\pm)$-diagonalizable if
$e$ acts by positive and negative eigenvalues on the objects of $\mathcal C$.
We define the zero of $\g$ to be a locally nilpotent, as well as a diagonalizable one-parameter element.\vspace*{1ex}
The following theorem is a particular case of Theorem 2.39 of \cite{M4}: 

\begin{Theorem}\label{gone} Let the Lie algebra $\g$ be generated by a set
  $E$ of one-parameter elements. Then $\mathcal C$, ${\mathcal C}^{du}$ is
  very good for integrating $\g$. Moreover, already the subgroup
\begin{eqnarray*}
  G_E := \bigcup_{m\in\N}\quad \bigcup_{e_1,\,e_2,\,\ldots,\,e_m\in E\setminus\{0\}} V_{e_1}V_{e_2}\cdots V_{e_m}
\end{eqnarray*} 
of $M^\times$ is dense in $M$. The coordinate ring $\FK{M}$ of $M$ is isomorphic to
its restriction $\FK{G_E}$ onto $G_E$ by the restriction map. 
\end{Theorem}
\begin{Remark}\label{Gfree} The group $G_E$ can be constructed in the same way as the Kac-Moody group in \cite{KP1}, or as the 
group associated to an integrable Lie algebra in \cite{K1}, \S 1.5, and \S 1.8:
Define the free product of groups $G_{free}:=\bigstar_{e\in E\setminus\{0\}} \F_e $, where $\F_e:=(\F,+)$ if $e$ is a locally nilpotent one-parameter element, and 
$\F_e:=(\F^\times,\,\cdot\,)$ if $e$ is a diagonalizable one-parameter element. Since the groups $V_e$, $e\in E\setminus\{0\}$, act on every $\g$-module
contained in $\mathcal C$, the homomorphisms $\kappa_e: \F_e\to V_e$, $e\in E\setminus\{0\}$, induce an action of $G_{free}$ on these $\g$-modules.
The homomorphisms $\kappa_e: \F_e\to V_e\subseteq G_E$, $e\in E\setminus\{0\}$, also induce a surjective homomorphism of groups 
$\kappa_{free}:G_{free}\to G_E$. Therefore $G_E$ is isomorphic to $G_{free}$
factored by the kernel of $\kappa_{free}$, which consists of the elements of $G_{free}$ which act trivially on every $\g$-module 
contained in $\mathcal C$.
\end{Remark}

To define, and in particular to work with the algebra of regular functions on $G_E$ we first introduce some notation. 
For every $e\in E\setminus\{0\}$ equip $V_e$ with its coordinate ring $\FK{V_e}$ obtained by restricting the coordinate ring
$\FK{M}$ of $M$ onto $V_e$. For every $p\in\N$ and $\underline{e}=(e_1,\,\ldots,\,e_p)\in (E\setminus\{0\})^p$ equip
\begin{eqnarray*}
 V_{\underline{e}}:=V_{e_1}\times\cdots\times V_{e_p} 
\end{eqnarray*}
with its coordinate ring $\FK{V_{\underline{e}}}$ as a product, given by $\FK{V_{e_1}}\otimes\cdots\otimes\FK{V_{e_p}}$. Let 
\begin{eqnarray*}
  m_{\underline{e}}: V_{\underline{e}}\to G_E .
\end{eqnarray*}
be the map given by multiplication. The image of $m_{\underline{e}}$ is $V_{e_1}V_{e_2}\cdots V_{e_p}$.\vspace*{1ex}

If $e\in E\setminus\{0\}$ is a locally nilpotent one-parameter element set $\F_e:=\F$ and equip $\F_e$ with the algebra $\FK{\F_e}:=\FK{t}$ of
polynomial functions. Set $\Z_e:=\Nn$. 

If $e\in E\setminus\{0\}$ is a diagonalizable one-parameter element set $\F_e:=\F^\times$. Recall that $\Z_e$ is the submonoid of $\Z$ 
given by the eigenvalues of $e_V$ for all objects $V$ of $\mathcal C$. Equip $\F_e$ with the algebra of functions
\begin{eqnarray*}
   \FK{\F_e}:=\bigoplus_{b\in \Z_e}\F \,t^b\;.
\end{eqnarray*}
For every $p\in\N$ and $\underline{e}=(e_1,\,\ldots,\,e_p)\in (E\setminus\{0\})^p$ equip
\begin{eqnarray*}
 \F_{\underline{e}}:=\F_{e_1}\times\cdots\times\F_{e_p} 
\end{eqnarray*}
with its algebra of functions $\FK{\F_{\underline{e}}}$ as a product, given by $\FK{\F_{e_1}}\otimes\cdots\otimes\FK{\F_{e_p}}$.
The elements of $\FK{\F_{\underline{e}}}$ are linear
combinations of monomials $t_1^{k_1}t_2^{k_2}\cdots t_p^{k_p}$ where $k_1\in\Z_{e_1},\,k_2\in\Z_{e_2},\,\ldots,\,k_p\in\Z_{e_p}$.
Define the map 
\begin{eqnarray*}
  \kappa_{\underline{e}}:\F_{\underline{e}}\to G_E\quad\mb{ by }\quad
  \kappa_{\underline{e}}(\underline{t}):=\kappa_{e_1}(t_1)\kappa_{e_2}(t_2)\cdots\kappa_{e_p}(t_p)
  \quad\mb{ where}\quad \underline{t}:=(t_1,\,\ldots,\,t_p)\in\F_{\underline{e}}.
\end{eqnarray*}
The image of $\kappa_{\underline{e}}$ is $V_{e_1}V_{e_2}\cdots V_{e_p}$.

\begin{Definition}\label{DefR} Let the Lie algebra $\g$ be generated by a set $E$ of one-parameter elements.
A function $f:G_E\to\F$ is called regular if for all $p\in\N$, and $\underline{e}\in (E\setminus\{0\})^p$ the following equivalent conditions are satisfied:
\begin{itemize}
\item[(i)]  $m_{\underline{e}}^*(f):=f\circ m_{\underline{e}}\in\FK{V_{\underline{e}}}$.
\item[(ii)] $\kappa_{\underline{e}}^*(f):=f\circ \kappa_{\underline{e}}\in\FK{\F_{\underline{e}}}$.
\end{itemize}
Denote by $\FK{G_E}_r$ the set of regular functions.
\end{Definition}

The next Proposition is proved in a similar way as the corresponding results in \cite{KP2} for Kac-Moody groups, and in \cite{K1}, 
\S 1.8 for the groups associated to linear integrable Lie algebras.
\begin{Proposition}\label{PropR} (a) $\FK{G_E}_r$ is a coordinate ring without zero divisors. Left and right multiplications with elements 
of $G_E$ induce comorphisms of $\FK{G_E}_r$.  If every non-zero diagonalizable one-parameter element of $E$ is $(\pm)$-diagonalizable, 
then the inverse map induces a comorphism of $\FK{G_E}_r$.

(b) The coordinate ring $\FK{G_E}$ is a subalgebra of $\FK{G_E}_r$.
\end{Proposition}

\begin{Remarks} (1) For the examples we treat in this
paper we reach a characterization of $\FK{G_E}$ as a subalgebra of $\FK{G_E}_r$ only by using the action of $G_E^{op}\times G_E$ on $\FK{G_E}_r$.

(2)  The coordinate ring $\FK{G_E}_r$ is the smallest coordinate ring on $G_E$, such that for all $p\in\N$,
and $e_1,\,\ldots,\,e_p\in E\setminus\{0\}$ the map
\begin{eqnarray*}
   m_{(e_1,\,\ldots,\,e_p)}:\;V_{e_1}\times\cdots\times V_{e_p}\to G_E
\end{eqnarray*}
given by multiplication is a morphism of sets with coordinate rings. 
\end{Remarks}

Next we introduce a subalgebra $\FK{U(\g)}_r$ of $U(\g)^*$, the algebra of regular linear functions on $U(\g)$, which contains the algebra of 
matrix coefficients $\FK{U(\g)}$ on $U(\g)$.
The definition of the linear regular functions on $U(\g)$ is the differentiated counterpart of the
definition of the regular functions on $G_E$. We first introduce the corresponding structures and notations:

Let $e\in E\setminus\{0\}$. Then $\F e$ is an abelian subalgebra of $\g$. Its universal enveloping algebra $U(\F e)$, and the dual $U(\F e)^*$ can be
described as follows:
\begin{eqnarray*}
  U(\F e)=\bigoplus_{n\in\Nn} \F \,\frac{e^n}{n!}\quad\mb{ and }\quad U(\F
  e)^*=\prod_{n\in\Nn} \F \tau^n\quad\mb{ where }\quad \tau^m(\frac{e^n}{n!})=\delta_{m n}.
\end{eqnarray*}
Since $U(\F e)$ is abelian the $U(\F e)$-action on $U(\F e)^*$ coincides with the
$U(\F e)^{op}$-action. We denote this action by $\mb{}_\diamond$. 
The element $e$ acts by formal differentiation, i.e.,
\begin{eqnarray}\label{xtau}
  e_\diamond\sum_{n\in\Nn}c_n \tau^n=\sum_{n\in\N} n c_n \tau^{n-1} \quad\mb{ where }\quad c_n\in\F,\;n\in\Nn.
\end{eqnarray}

Let $e\in E\setminus\{0\}$ be a locally nilpotent one-parameter element. As we have seen in the proof of Theorem 2.35 of \cite{M4}, the algebra of matrix coefficients on 
$U(\F e)$ is given by
\begin{eqnarray*}
   \FK{U(\F e)}=\bigoplus_{n\in\Z_e}\F \tau^n\subseteq U(\F e)^*.
\end{eqnarray*}
(Recall that $\Z_e=\Nn$). Since $e$ acts locally nilpotent on $\FK{U(\F  e)}$, also the group $V_e$ acts on $\FK{U(\F  e)}$. From equation (\ref{xtau}) follows that 
this action is given on the  elements $\tau^n\in U(\g)^*$, $n\in\Z_e$, by
\begin{eqnarray*}
   \exp(te)_\diamond \tau^n =(t  +\tau)^n\quad \mb{ where }\quad t\in\F,\;n\in\Nn.
\end{eqnarray*}

Let $e\in E\setminus\{0\}$ be a diagonalizable one-parameter element. As we have seen in the proof of Theorem 2.37 of \cite{M4}, the
algebra of matrix coefficients on $U(\F e)$ is given by 
\begin{eqnarray*}
   \FK{U(\F e)}=\bigoplus_{n\in\Z_e}\F \exp(\tau)^n\subseteq U(\F e)^*\quad\mb{
   where }\quad \exp(\tau)^n=\exp(n\tau):=\sum_{k\in\Nn}\frac{n^k}{k!} \tau^k\in U(\g)^*.
\end{eqnarray*}
(Recall that $\Z_e$ is the submonoid of $(\Z,+)$ given by the eigenvalues of $e_V$ for all objects $V$ of $\mathcal C$.) Since $e$ acts 
diagonalizable on on $\FK{U(\F  e)}$, also the group $V_e$ acts on $\FK{U(\F e)}$. From equation (\ref{xtau}) follows that this action is 
given on the  elements  $\exp(\tau)^n\in U(\g)^*$, $n\in\Z$, by
\begin{eqnarray*}
   (t^e)_\diamond \exp(\tau)^n = \left(t\exp(\tau)\right)^n\quad \mb{ where }\quad t\in\F^\times,\;n\in\Z.
\end{eqnarray*}

To cut short our notation we define
\begin{eqnarray}\label{etadef}
\eta:=\left\{\begin{array}{cl}
    \tau &\mb{ if } e \mb{ is a locally nilpotent one-parameter element }\\
    \exp(\tau) &\mb{ if } e \mb{ is a diagonalizable one-parameter element }
             \end{array}\right\}\in \FK{U(\F e)}.
\end{eqnarray}

The universal enveloping algebras $U(\F e)$, $e\in E\setminus\{0\}$, embed into $U(\g)$ in the obvious way. Let $p\in\N$ and $\underline{e}=(e_1,\,\ldots,\,e_p)\in
(E\setminus\{0\})^p$. Let
\begin{eqnarray*}
  \rho_{\underline{e}}:\,U(\F e_1)\otimes \cdots\otimes U(\F e_p)\to U(\g)
\end{eqnarray*}
be the linear map given by $\rho_{\underline{e}}(y_1\otimes\cdots\otimes y_p):=y_1\cdots y_p$,
where $y_1\in U(\F e_1)$, \ldots, $y_p\in U(\F e_p)$. For $e\in E\setminus\{0\}$ set $\rho_e:=\rho_{(e)}$.

Note that the space of linear functions on $U(\F e_1)\otimes \cdots\otimes U(\F e_p)$ can be described by
\begin{eqnarray*}
  \left(U(\F e_1)\otimes \cdots\otimes U(\F e_p)\right)^*=\prod_{k_1,\,\ldots,\, k_n\in\Nn}\F \,\tau_1^{k_1}\otimes\cdots\otimes\tau_p^{k_p}.
\end{eqnarray*}
Identify $\FK{U(\F e_1)}\otimes \cdots\otimes \FK{U(\F e_p)}$ with the corresponding subspace of $\left(U(\F e_1)\otimes \cdots\otimes U(\F e_p)\right)^*$.

\begin{Definition} A linear function $h\in U(\g)^*$ of $U(\g)$ is called regular if for all $\underline{e}=(e_1,\,\ldots,\,e_p)\in (E\setminus\{0\})^p$, $p\in\N$, we have
\begin{eqnarray*}
  h\circ\rho_{\underline{e}}\in \FK{U(\F e_1)}\otimes \cdots\otimes \FK{U(\F e_p)}.
\end{eqnarray*}
Denote by $\FK{U(\g)}_r$ the set of regular linear functions of $U(\g)$.
\end{Definition}

For $e_1,\,\ldots,\,e_p\in E\setminus\{0\}$, $p\in\N$, equip $\FK{U(\F e_1)}\otimes \cdots\otimes \FK{U(\F e_p)}$ with the algebra structure induced by $\FK{U(\F e_j)}$, 
$j=1,\,\ldots,\, p$. 

\begin{Proposition} \label{ureg}(a) $\FK{U(\g)}_r$ is a $U(\g)^{op}\otimes  U(\g)$-invariant subalgebra of $U(\g)^*$. 
If every non-zero diagonalizable one-parameter element of $E$ is $(\pm)$-diagonalizable, then $\FK{U(\g)}_r$ is also invariant under the 
dual of the antipode of $U(\g)$. 

(b) For $\underline{e}=(e_1,\,\ldots,\,e_p)\in (E\setminus\{0\})^p$, $p\in\N$, the map
\begin{eqnarray*}
   \rho_{\underline{e}}^*:\FK{U(\g)}_r \to \FK{U(\F e_1)}\otimes \cdots\otimes
   \FK{U(\F e_p)} 
\end{eqnarray*}
defined by $\rho_{\underline{e}}^*(h):=h\circ\rho_{\underline{e}}$, where $h\in\FK{U(\g)}_r$, is a morphism of algebras.

(c) The algebra of matrix coefficients $\FK{U(\g)}$ on $U(\g)$ is a subalgebra of $\FK{U(\g)}_r$.  
\end{Proposition}

\Proof To (a) and (b): It is easy to check that $\FK{U(\g)}_r$ is a $U(\g)^{op}\otimes U(\g)$-invariant linear subspace of $U(\g)^*$. 
Denote by $\cdot$ the multiplication of $U(\g)^*$. To show that $\FK{U(\g)}_r$ is a subalgebra of $U(\g)^*$, and to prove part (b) of the 
proposition, it is sufficient to show: Let $p\in\N$ and $\underline{e}=(e_1,\,\ldots,\,e_p)\in (E\setminus\{0\})^p$. Then  
\begin{eqnarray*}
  1\circ  \rho_{\underline{e}}=1\otimes\cdots\otimes 1 \quad\mb{ and }\quad
  (h\cdot\ti{h})\circ\rho_{\underline{e}}= (h\circ\rho_{\underline{e}})\bullet(\ti{h}\circ\rho_{\underline{e}}),
\end{eqnarray*}
where $h,\ti{h}\in \FK{U(\g)}_r$, and $\bullet$ denotes the multiplication of $\FK{U(\F e_1)}\otimes\cdots\otimes \FK{U(\F e_p)}$

The first equation is trivial. For $\underline{k}:=(k_1,\,\ldots,\, k_p)\in\Nn^p$ define
$e_{\underline{k}}:=\frac{e_1^{k_1}\cdots e_p^{k_p}}{k_1!\cdots k_p!}$. These
elements satisfy
\begin{eqnarray*}
   \Delta(e_{\underline{k}})=\sum_{\underline{a},\underline{b}\in\Nn^p,\;\underline{a}+\underline{b}= \underline{k}}
   e_{\underline{a}}\otimes e_{\underline{b}}.
\end{eqnarray*}
For $\underline{k}:=(k_1,\,\ldots,\, k_p)\in\Nn^p$ define $\tau^{\underline{k}}:=\tau_1^{k_1}\otimes\cdots\otimes
\tau_p^{k_p}$. We have
\begin{eqnarray*}
  (h\cdot \ti{h})\circ \rho_{\underline{e}}=\sum_{\underline{k}\in\Nn^p}\left((h\cdot\ti{h})(e_{\underline{k}})\right)\tau^{\underline{k}}
  =\sum_{\underline{k}\in\Nn^p}\left((h\otimes\ti{h})(\Delta(e_{\underline{k}}))\right)\tau^{\underline{k}}\\
  = \sum_{\underline{k}\in\Nn^p}\;\sum_{\underline{a},\underline{b}\in\Nn^p,\;\underline{a}+\underline{b}=\underline{k}}
     h(e_{\underline{a}})\ti{h}(e_{\underline{b}})\tau^{\underline{k}}
  = \sum_{\underline{k}\in\Nn^p}\;\sum_{\underline{a},\underline{b}\in\Nn^p,\;\underline{a}+\underline{b}= \underline{k}}
      h(e_{\underline{a}})\tau^{\underline{a}}\bullet\ti{h}(e_{\underline{b}})\tau^{\underline{b}}\\
  = \left(\sum_{\underline{a}\in\Nn^p} h(e_{\underline{a}})\tau^{\underline{a}}\right)\bullet 
               \left(\sum_{\underline{b}\in\Nn^p}\ti{h}(e_{\underline{b}})\tau^{\underline{b}}\right)=
       (h\circ \rho_{\underline{e}})\bullet(\ti{h}\circ \rho_{\underline{e}}).
\end{eqnarray*}

Now assume that every non-zero diagonalizable one-parameter element of $E$ is
$(\pm)$-diagonalizable. Denote by $S$ be the antipode of $U(\g)$. Let $h\in \FK{U(\g)}_r$. Let $p\in\N$ and $e_1,\,\ldots,\,e_p\in
E\setminus\{0\}$. Recall the notation (\ref{etadef}). Since $h$ is regular,
\begin{eqnarray*}
    h\circ\rho_{(e_p,\,\ldots,\,e_1)}=\sum_{finitely\; many\atop k_1,\,\ldots,\,k_p}c_{k_p\cdots
    k_1} \eta_p^{k_p}\otimes\cdots\otimes \eta_1^{k_1}\in \FK{U(\F e_p)}\otimes \cdots\otimes \FK{U(\F e_1)}.
\end{eqnarray*} 
Denote by $S_j$ the antipode of $U(\F e_j)$, $j=1,\,\ldots,\,p$. Then
\begin{eqnarray*}
   S^*(h)\circ \rho_{(e_1,\,\ldots,\,e_p)}=\sum_{finitely \;many\atop k_1,\,\ldots,\,k_p} c_{k_p\cdots k_1} 
   S_1^*(\eta_1^{k_1})\otimes\cdots\otimes S_p^*(\eta_p^{k_p}).
\end{eqnarray*}
Because of $S_j^*(\tau_j)=-\tau_j$ and $S_j^*(\exp(\tau_j))=\exp(-\tau_j)$ the sum on the right is contained in 
$\FK{U(\F e_1)}\otimes \cdots\otimes \FK{U(\F e_p)}$.\vspace*{1ex}

To (c): Let $V$ be an object of $\mathcal C$. Let $v\in V$ and $\phi\in
V^{du}$. Let $p\in\N$ and $e_1,\,\ldots,\,e_p\in E\setminus\{0\}$. Let $y_1\in
U(\F e_1)$, \ldots, $y_p\in U(\F e_p)$. 

If $e_p$ acts locally nilpotent there exists an integer $N_p\in\N$ such that
$e_p^{N_p + 1}v=0$. We get the following development in position $p$:
\begin{eqnarray}\label{dp1}
 \phi(y_1\cdots y_{p-1}y_p v)=\sum_{k_p=0}^{N_p}\phi(y_1\cdots y_{p-1}e_p^{k_p}v)\,\frac{\tau_p^{k_p}}{k_p!}(y_p).
\end{eqnarray}
If $e_p$ acts diagonalizable there exists an integer $N_p\in\Nn$ such that
$v=\sum_{n\in \Z_e,\,|n|\leq N_p} v_n$, where $v_n$ is either an eigenvector of $e$ on $V$ to the
eigenvalue $n$ or the zero vector. We get the following development in position $p$:
\begin{eqnarray}\label{dp2}
 \phi(y_1\cdots y_{p-1}y_p v)=\sum_{k_p\in \Z_e,\,|k_p|\leq N_p}\phi(y_1\cdots y_{p-1}v_{k_p})\left(\exp(k_p\tau_p)(y_p)\right).
\end{eqnarray}

Developing (\ref{dp1}) resp. (\ref{dp2}) successively in the positions $p-1$ to $1$ we find $g_{\phi v}\circ \rho_{(e_1,\,\ldots,\,e_p)}\in
\FK{U(\F e_1)}\otimes\cdots\otimes \FK{U(\F e_p)}$.
\qed

At first the algebra $\FK{U(\g)}_r$ looks even more unpleasant than $\FK{G_E}_r$. But actually it is easier to work with. In particular, its
definition boils down nicely in the important case where $E$ consists only of locally nilpotent one-parameter elements:
\begin{Proposition}\label{nilregalg} Suppose that $E$ consists of locally nilpotent one-parameter elements. 
Then $\FK{U(\g)}_r$ is the subalgebra of $U(\g)^*$ given by the elements $h\in U(\g)^*$, such that for all $p\in\N$, and for all $e_1,\,\ldots,\,e_p\in E\setminus\{0\}$ we have
\begin{eqnarray*}
  h(e_1^{k_1}\cdots e_p^{k_p})\neq 0 \quad\mb{ for at most finitely many }k_1,\,\ldots,\,k_p\in\Nn.
\end{eqnarray*}
\end{Proposition}

\Proof Let $h\in U(\g)^*$. For $\underline{e}=(e_1,\,\ldots,\,e_p)\in (E\setminus\{0\})^p$ and $p\in\N$ we can develop
\begin{eqnarray*}
  h\circ \rho_{\underline{e}} =
    \sum_{k_1,\,\ldots,\,k_p\in\Nn} h(\frac{e_1^{k_1}}{k_1!}\cdots\frac{e_p^{k_p}}{k_p!})\,\tau_1^{k_1}\otimes\cdots\otimes\tau_p^{k_p}
    \in \left(U(\F e_1)\otimes\cdots\otimes U(\F e_p)\right)^*.
\end{eqnarray*}
The sum on the right is contained in $\FK{U(\F e_1)}\otimes \cdots\otimes \FK{U(\F e_p)}$ if and only if
$h(\frac{e_1^{k_1}}{k_1!}\cdots\frac{e_p^{k_p}}{k_p!} )\neq 0$ for at most finitely many $k_1,\,\ldots,\, k_p\in\Nn$.
\qed

Next we investigate when the counterparts $\FK{G_E}_r$ and $\FK{U(\g)}_r$ are isomorphic. If $\FK{G_E}_r$ and $\FK{U(\g)}_r$ are isomorphic algebras, 
then there exists a natural action of $U(\g^{op})\otimes U(\g)$ on $\FK{G_E}_r$, as well as a natural action of $G_E^{op}\times G_E$ on $\FK{U(\g)}_r$. 
In Theorem \ref{RegDual} we show that this is already sufficient. As a preparation we first show Proposition \ref{eact} and Theorem
\ref{Vact}, which are quite useful.

For $e\in E\setminus\{0\}$ define
\begin{eqnarray*}
  \frac{d}{dt}\res{e} := \left\{\begin{array}{cl}\frac{d}{dt}\res{t=0} &\mb{ if } e \mb{ is a locally nilpotent one-parameter element}\\
    \frac{d}{dt}\res{t=1}  &\mb{ if } e \mb{ is a diagonalizable one-parameter element}
        \end{array}\right. .
\end{eqnarray*} 
Recall the definition of the element $\eta$ from (\ref{etadef}).

\begin{Proposition}\label{eact} (a) For every $e\in E\setminus\{0\}$ we get a left invariant derivation of $\FK{G_E}_r$ by  
\begin{eqnarray}\label{eaction}
    (e_\ro f)(g):=\frac{d}{dt}\res{e}f(g\kappa_e(t)) \quad\mb{ where } \quad f\in\FK{G_E}_r, \;g\in G_E.
\end{eqnarray}
Similarly, we get a right invariant derivation of $\FK{G_E}_r$ by  
\begin{eqnarray}\label{eaction2}
    (e_\lo f)(g):=\frac{d}{dt}\res{e}f(\kappa_e(t)g) \quad\mb{ where } \quad f\in\FK{G_E}_r, \;g\in G_E.
\end{eqnarray} 
For $e,\ti{e}\in E\setminus\{0\}$ the derivations $e_\ro$ and $\ti{e}_\lo$ commute.\vspace*{1ex}

(b)  Let $f\in \FK{G_E}_r$. Let $e_1,\,\ldots,\,e_q \in E\setminus\{0\}$, $q\in\N$. Develop
\begin{eqnarray}\label{eaction3}
   f(\kappa_{e_1}(t_1)\cdots\kappa_{e_q}(t_q))=\sum_{n_1,\,\ldots,\,n_q}
     c_{n_1 \cdots n_q } t_1^{n_1}\cdots t_q^{n_q} \quad\mb{
     where }\quad t_1\in\F_{e_1},\,\cdots,\,t_q\in\F_{e_q}.
\end{eqnarray}
Let $p\in\{1,\,2,\,\ldots,\,q\}$. Then for all $t_1\in\F_{e_1},\,\cdots,\,t_{p-1}\in\F_{e_{p-1}}$, $k_p,\,\ldots,\,k_q\in\Nn$
we have 
\begin{eqnarray}
  ({e_p^{k_p}}_\ro (\cdots ({e_q^{k_q}}_\ro  f)\cdots))(\kappa_{e_1}(t_1)\cdots\kappa_{e_{p-1}}(t_{p-1}))\qquad\qquad\qquad\qquad\quad\nonumber\\
    = \sum_{n_1,\,\ldots,\,n_q} c_{n_1\cdots
 n_q}\,\eta_p^{n_p}(e_p^{k_p}) \cdots \eta_q^{n_q}(e_q^{k_q})\,t_1^{n_1}\cdots
 t_{p-1}^{n_{p-1}}.\label{eaction4a}
\end{eqnarray}
(For $p=1$ replace $\kappa_{e_1}(t_1)\cdots\kappa_{e_{p-1}}(t_{p-1})$ by $1$,
and omit the expression $t_1^{n_1}\cdots t_{p-1}^{n_{p-1}}$.)
For all $t_{p+1}\in\F_{e_{p+1}},\,\cdots,\,t_q\in\F_{e_q}$, $k_1,\,\ldots,\,k_p\in\Nn$, we have
\begin{eqnarray}
({e_p^{k_p}}_\lo (\cdots ({e_1^{k_1}}_\lo f)\cdots))(\kappa_{e_{p+1}}(t_{p+1})\cdots\kappa_{e_q}(t_q))\qquad\qquad\qquad\qquad\quad \nonumber\\
= \sum_{n_1,\,\ldots,\,n_q} c_{n_1\cdots n_q}\,\eta_1^{n_1}(e_1^{k_1}) \cdots \eta_p^{n_q}(e_p^{k_p})\,t_{p+1}^{n_{p+1}}\cdots t_q^{n_q}\label{eaction4b}.
\end{eqnarray}
(For $p=q$ replace $\kappa_{e_{p+1}}(t_{p+1})\cdots\kappa_{e_q}(t_q)$ by $1$, and omit the expression $t_{p+1}^{n_{p+1}}\cdots t_q^{n_q}$.)
\end{Proposition}

\Proof We first show part (b) of the theorem. Let $e_p$ be a locally nilpotent one-parameter element. By
(\ref{eaction3}) we find
\begin{eqnarray*}
 \left({e_p^{k_p}}_\ro f \right)\left(\kappa_{(e_1,\,\ldots,\,e_{p-1})}(t_1,\,\ldots,\, t_{p-1})\right)\\
  =\frac{d}{d u_1}\res{u_1=0}\cdots \frac{d}{d u_{k_p}}\res{u_{k_p}=0}
          f(\kappa_{(e_1,\,\ldots,\,e_{p-1})}(t_1,\,\ldots,\, t_{p-1})\exp(u_1 e_p)\cdots \exp(u_{k_p} e_p))\\
  =  \frac{d}{d u_1}\res{u_1=0}\cdots \frac{d}{d u_{k_p}}\res{u_{k_p}=0}\sum_{n_1\cdots n_p}c_{n_1\cdots
          n_p}t_1^{n_1}\cdots t_{p-1}^{n_{p-1}}(u_1+\cdots +u_{k_p})^{n_p}\\
  =  \sum_{n_1\cdots n_p}c_{n_1\cdots n_p}t_1^{n_1}\cdots t_{p-1}^{n_{p-1}} \delta_{n_p k_p} n_p! 
  =  \sum_{n_1\cdots n_p}c_{n_1\cdots n_p}t_1^{n_1}\cdots t_{p-1}^{n_{p-1}}\tau^{n_p}(e_p^{k_p}).
\end{eqnarray*}
for all $k_p\in\N$. Also the first and the last expression of this chain of equalities coincide for $k_p=0$.
If $e_p\in E\setminus\{0\}$ is a diagonalizable one-parameter element, then by
a similar computation
\begin{eqnarray*}
 \left({e_p^{k_p}}_\ro f \right)\left(\kappa_{(e_1,\,\ldots,\,e_{p-1})}(t_1,\,\ldots,\, t_{p-1})\right)
     %\\=\frac{d}{d u_1}\res{u_1=1}\cdots \frac{d}{d u_{k_p}}\res{u_{k_p}=1}
     %        f(\kappa_{(e_1,\,\ldots,\,e_{p-1})}(t_1,\,\ldots,\, t_{p-1}) u_1^{e_p} \cdots u_{k_p}^{e_p})\\
     %=  \frac{d}{d u_1}\res{u_1=1}\cdots \frac{d}{d u_{k_p}}\res{u_{k_p}=1}\sum_{n_1\cdots n_p}c_{n_1\cdots
     %        n_p}t_1^{n_1}\cdots t_{p-1}^{n_{p-1}}(u_1\cdots u_{k_p})^{n_p}\\
     %=  \sum_{n_1\cdots n_p}c_{n_1\cdots n_p}t_1^{n_1}\cdots t_{p-1}^{n_{p-1}} n_p^{k_p}  
=  \sum_{n_1\cdots n_p}c_{n_1\cdots n_p}t_1^{n_1}\cdots t_{p-1}^{n_{p-1}}\left(\exp(n_p\tau)(e_p^{k_p})\right)
\end{eqnarray*}
for all $k_p\in\Nn$.
%Also the first and the last expression of this chain of equalities coincide for $k_p=0$. 
Applying these steps successively $(q-p)$-times 
we get formula (\ref{eaction4a}). Formula (\ref{eaction4b}) is proved in a similar way.\vspace*{1ex}

Now part (a) of the theorem follows from part (b) and from standard properties of the differentiation.
\qed

\begin{Theorem}\label{Vact} (a) For every $e\in E\setminus \{0\}$ the group $V_e$ acts by morphisms
of algebras on $\FK{U(\g)}_r$ by
\begin{eqnarray}\label{Vaction}
   (g_\ro h)(y):=\left(g_\diamond (\rho_e^*(y_\lo h))\right)(1) \quad\mb{
   where }\quad g\in V_e,\;h\in \FK{U(\g)}_r,\;y\in U(\g).
\end{eqnarray}
For every $e\in E\setminus \{0\}$ the group $V_e^{op}$ acts by morphisms of algebras on $\FK{U(\g)}_r$ by
\begin{eqnarray}\label{Vaction2}
   (g_\lo h)(y):=\left(g_\diamond (\rho_e^*(y_\ro h))\right)(1) \quad\mb{ where }\quad g\in V_e,\;h\in \FK{U(\g)}_r,\;y\in U(\g).
\end{eqnarray}
For $e,\,\ti{e}\in E\setminus \{0\}$ the actions of $V_e$ and $V_{\ti{e}}^{op}$ commute. \vspace*{1ex}

(b) Let $h\in \FK{U(\g)}_r$. Let
$e_1,\,\ldots,\,e_q \in E\setminus\{0\}$, $q\in\N$. Develop
\begin{eqnarray}\label{Vaction3}
  \rho_{(e_1,\,\cdots,\, e_q)}^*(h)=\sum_{k_1,\,\ldots,\,k_q}  c_{k_1 \cdots k_q }\eta_1^{k_1}\otimes\cdots\otimes\eta_q^{k_q}.
\end{eqnarray}
Let $p\in\{1,\,2,\,\ldots,\,q\}$. Then for all $t_p\in\F_{e_p},\,\cdots,\,t_q\in\F_{e_q}$ we have
\begin{eqnarray}\label{Vaction4a}
\rho_{(e_1,\,\cdots,\, e_{p-1})}^*\left(\, \kappa_{e_p}(t_p)_\ro(\cdots(\kappa_{e_q}(t_q)_\ro h )\cdots)\,\right)=\sum_{k_1,\,\ldots,\,k_q}
     c_{k_1 \cdots k_q }  t_p^{k_p}\cdots
     t_q^{k_q}\,\eta_1^{k_1}\otimes\cdots\otimes\eta_{p-1}^{k_{p-1}}.
\end{eqnarray}
(For $p=1$ replace $\rho_{(e_1,\,\cdots,\, e_{p-1})}^*$ by the evaluation map in $1$, and omit  
$\eta_1^{k_1}\otimes\cdots\otimes\eta_{p-1}^{k_{p-1}}$.) For all $t_1\in\F_{e_1},\,\cdots,\,t_p\in\F_{e_p}$ we have
\begin{eqnarray}
  \rho_{(e_{p+1},\,\cdots,\, e_q)}^*\left(\,\kappa_{e_p}(t_p)_\lo(\cdots(\kappa_{e_1}(t_1)_\lo h )\cdots)\,\right)=\sum_{k_1,\,\ldots,\,k_q}
 c_{k_1 \cdots k_q} t_1^{k_1}\cdots t_p^{k_p}\,\eta_{p+1}^{k_{p+1}}\otimes\cdots\otimes \eta_q^{k_q}.\label{Vaction4b}
\end{eqnarray}
(For $p=q$ replace $\rho_{(e_{p+1},\,\cdots,\, e_q)}^*$ by the evaluation map in $1$, and omit 
$\eta_{p+1}^{k_{p+1}}\otimes\cdots\otimes \eta_q^{k_q}$.)
\end{Theorem}

\Proof We first show part (b) of the theorem. By (\ref{Vaction3}) we find
\begin{eqnarray*}
   ((y_1\cdots y_{q-1})_\lo h)\circ\rho_{e_q}=\sum_{k_1,\,\ldots,\,k_q}c_{k_1\cdots
   k_q}\eta_1^{k_1}(y_1)\cdots \eta_1^{k_{q-1}}(y_{q-1}) \, \eta_q^{k_q}.
\end{eqnarray*}
for all $y_1\in U(\F e_1)$, \ldots, $y_{q-1}\in U(\F e_{q-1})$. It follows
\begin{eqnarray*}
 (\kappa_{e_q} (t_q)_\ro h)(y_1\cdots y_{q-1})
   =\left(\kappa_{e_q}(t_q)_\diamond\left(( (y_1\cdots y_{q-1})_\lo h)\circ \rho_{e_q} \right)\right)(1)\\
   =\sum_{k_1,\,\ldots,\,k_q}c_{k_1\cdots k_q}\eta_1^{k_1}(y_1)\cdots \eta_{q-1}^{k_{q-1}}(y_{q-1})
   \underbrace{(\kappa_{e_q}(t_q)_\diamond \eta^{k_q})(1)}_{=t_q^{k_q}}
\end{eqnarray*}
for all $y_1\in U(\F e_1)$, \ldots, $y_{q-1}\in U(\F e_{q-1})$. Equivalently,
\begin{eqnarray*}
  (\kappa_{e_q} (t_q)_\ro h)\circ \rho_{(e_1,\,\ldots,\,e_{q-1})}
   =\sum_{k_1,\,\ldots,\,k_q}c_{k_1\cdots k_p}t_q^{k_q}\,\eta_1^{k_1}\otimes\cdots \otimes \eta_{q-1}^{k_{q-1}} 
\end{eqnarray*}
Repeating the same step successively $(q-p)$-times gives formula (\ref{Vaction4a}). Formula (\ref{Vaction4b}) is proved in a similar way.\vspace*{1ex}

Now we can prove part (a) of the theorem. If $g\in V_e$ and $h\in \FK{U(\g)}_r$ then also $g_\ro h\in \FK{U(\g)}_r$, resp. $\,g_\lo h\in \FK{U(\g)}_r$. 
This follows from formula (\ref{Vaction4a}) for $e_q=e$ and $p=q-1$, resp. (\ref{Vaction4b}) for $e_1=e$ and $p=1$.\vspace*{1ex}

We only show that $V_e$ acts on $\FK{U(\g)}_r$, the proof for $V_e^{op}$ is similar. Let $y\in U(\g)$ and $h\in \FK{U(\g)}_r$. Develop
\begin{eqnarray}\label{developact}
         (y_\lo h)\circ \rho_{(e,e)}=\sum_{m_1,\, m_2} \ti{c}_{m_1 m_2} \eta^{m_1}\otimes \eta^{m_2} 
                \quad\mb{ and }\quad      
         (y_\lo h)\circ \rho_e=\sum_n c_n \eta^n .
\end{eqnarray}
By part (b), formulas (\ref{Vaction4a}) and (\ref{Vaction4b}), we find 
\begin{eqnarray}\label{doublekappa}
 \left(\kappa_e(t)_\ro (  \kappa_e(\ti{t})_\ro h)\right)(y)=\sum_{m_1 m_2}\ti{c}_{m_1 m_2}t^{m_1}\ti{t}^{m_2}.
\end{eqnarray}

First let $e$ be a locally nilpotent one-parameter element. Let
$k_1,\,k_2\in\Nn$. Because of $(y_\lo h)(x^{k_1}x^{k_2})=(y_\lo h) (x^{k_1+k_2})$
we find from (\ref{developact}):
\begin{eqnarray*}
   k_1!\, k_2!\,\ti{c}_{k_1k_2} =\sum_{m_1,\, m_2}\ti{c}_{m_1
   m_2}\tau^{m_1}(x_1^{k_1})\tau^{m_2}(x_2^{k_2})=\sum_n c_n
  \tau^n(x^{k_1+k_2})=(k_1+k_2)! \,c_{k_1+k_2}.
\end{eqnarray*}
Inserting in (\ref{doublekappa}) we get
\begin{eqnarray*}
   \left(\kappa_e(t)_\ro (  \kappa_e(\ti{t})_\ro h)\right)(y)=\sum_{m_1,\,m_2\in\Nn} \ti{c}_{m_1 m_2} t^{m_1}\ti{t}^{m_2}
  = \sum_{m_1, \,m_2\in\Nn} c_{m_1+m_2}  \frac{(m_1+m_2)!}{m_1!m_2!}\,t^{m_1}\ti{t}^{m_2}\\
  =\sum_{n\in\Nn} c_n \sum_{m_1,\,m_2\in\Nn,\,m_1+m_2=n}\frac{(m_1+m_2)!}{m_1!m_2!}\,t^{m_1}\ti{t}^{m_2}  
  = \sum_{n\in\Nn}c_n (t+\ti{t})^n.
\end{eqnarray*}
for all $t,\,\ti{t}\in\F$. By part (b), formula (\ref{Vaction4a}), the last sum coincides with $(\kappa_e(t+\ti{t})_\ro h)(y)$.

Now let $e$ be a diagonalizable one-parameter element. Let
$k_1,\,k_2\in\Nn$. Because of $(y_\lo h)(x^{k_1}x^{k_2})=(y_\lo h)(x^{k_1+k_2})$
we find from (\ref{developact}):
\begin{eqnarray*}
     \sum_{m_1,\, m_2}\ti{c}_{m_1 m_2} m_1^{k_1} m_2^{k_2}  
   =\sum_{m_1,\, m_2}\ti{c}_{m_1 m_2}\left(\exp(m_1\tau)(x_1^{k_1})\right)\left(\exp(m_2\tau)(x_2^{k_2})\right)\\
   =\sum_n c_n\exp(n\tau)(x^{k_1+k_2})= \sum_n c_n n^{k_1+k_2}.
\end{eqnarray*}
We may write these equations in the form
\begin{eqnarray*}
   \sum_{m_1} \left(\sum_{m_2}(\ti{c}_{m_1 m_2}
   -\delta_{m_1 m_2}c_{m_1})m_2^{k_2}\right) m_1^{k_1}=0 
     \quad\mb{ for all }\quad k_1,\,k_2\in\Nn. 
\end{eqnarray*}
Since only finitely many of the coefficients $\ti{c}_{m_1 m_2}$ and $c_n$ are
non-zero, this is a interlocked double Van-der-Monde system. It follows 
\begin{eqnarray*}
  \ti{c}_{m_1 m_2}=0 \quad\mb{ for all} \quad m_1,\,m_2,\; m_1\neq m_2, \quad\mb{ and }\quad\ti{c}_{nn}=c_n\quad \mb{ for all } \quad n.
\end{eqnarray*}
Inserting in (\ref{doublekappa}) we find
\begin{eqnarray*}
   \left(\kappa_e(t)_\ro (  \kappa_e(\ti{t})_\ro h)\right)(y)=\sum_{m_1,\,m_2} \ti{c}_{m_1 m_2} t^{m_1}\ti{t}^{m_2}
  = \sum_n c_n t^n\ti{t}^n
\end{eqnarray*}
for all $t,\,\ti{t}\in\F^\times$. By part (b), formula (\ref{Vaction4a}), the last sum coincides with $(\kappa_e(t\ti{t})_\ro h)(y)$.\vspace*{1ex}

Let $h\in \FK{U(\g)}_r$. Fix $e,\,\ti{e}\in E\setminus \{0\}$. Let $e_1,\,\ldots,\,e_p\in E\setminus\{0\}$, $p\in\Nn$. Develop
\begin{eqnarray*}
 h\circ\rho_{(\ti{e},\,e_1,\,\ldots,\, e_p,\,e)}=\sum_{\ti{n},\,n_1,\,\ldots,\,n_p,\,n}c_{\ti{n} n_1\cdots n_p n}
 \ti{\eta}^{\ti{n}}\otimes \eta_1^{n_1}\otimes\cdots\otimes \eta_p^{n_p}\otimes \eta^n.
\end{eqnarray*}
By part (b), formulas (\ref{Vaction4a}) and (\ref{Vaction4b}), we find that 
$\left(\kappa_{\ti{e}}(\ti{t})_\lo (\kappa_e(t)_\ro  h)\right)\circ\rho_{(e_1,\,\ldots,\,e_p)}$ as well as 
$\left(\kappa_e (t)_\ro   (\kappa_{\ti{e}}(\ti{t})_\lo h)\right)\circ \rho_{(e_1,\,\ldots,\,e_p)}$ equal
\begin{eqnarray*}
  \sum_{\ti{n},\,n_1,\,\ldots,\,n_p,\,n}c_{\ti{n} n_1\cdots n_p n}\,\ti{t}^{\ti{n}}
  t^n\,\eta_1^{n_1}\otimes\cdots \otimes \eta_p^{n_p} \quad\mb{ for all }\quad\ti{t}\in\F_{\ti{e}},\;t\in\F_e.
\end{eqnarray*}
The Lie algebra $\g$ is generated by $E$. Therefore, the universal enveloping
algebra $U(\g)$ is generated by its subalgebras $U(\F e)$, $e\in E$. It follows
$\left(\kappa_{\ti{e}}(\ti{t})_\lo (\kappa_e(t)_\ro  h)\right)= 
\left(\kappa_e (t)_\ro   (\kappa_{\ti{e}}(\ti{t})_\lo h)\right)$ for all $\ti{t}\in\F_{\ti{e}}$, $t\in\F_e$.
\qed

\begin{Theorem}\label{RegDual} (a) Suppose there exists an action of the Lie algebras $\g$ on $\FK{G_E}_r$ extending the actions
(\ref{eaction}) of the elements $e\in E$, and an action of the Lie algebra $\g^{op}$ on $\FK{G_E}_r$ extending the actions (\ref{eaction2})
of the elements $e\in E$. Then we get a $U(\g)^{op}\otimes U(\g)$-equivariant embedding of algebras $\Phi:\FK{G_E}_r\to \FK{U(\g)}_r$ by
\begin{eqnarray}\label{phidef}
   \Phi(f)(x):=(x_\ro f)(1)=(x_\lo f)(1) \quad \mb{ where }\quad f\in\FK{G_E}_r,\;x\in U(\g).
\end{eqnarray}

(b) Suppose there exists an action of the group $G_E$ extending the actions (\ref{Vaction})  of the groups $V_e$, $e\in E$, 
and an action of $G_E^{op}$ on $\FK{U(\g)}_r$ extending the actions (\ref{Vaction2}) of the groups $V_e^{op}$, $e\in E$.
Then we get a $G_E^{op}\times G_E$-equivariant embedding of algebras $\Xi:\FK{U(\g)}_r \to \FK{G_E}_r$ by
\begin{eqnarray}\label{xidef}
   \Xi(h)(g):=(g_\ro h)(1)=(g_\lo h)(1) \quad \mb{ where }\quad h\in \FK{U(\g)}_r,\;g\in G_E.
\end{eqnarray}

(c) If the assumptions of (a) and (b) are satisfied, then $\Phi$ and $\Xi$ are inverse maps.
\end{Theorem}
\begin{Remark} (a) We have seen that $\FK{G_E}\subseteq \FK{G_E}_r$, and $\FK{U(\g)}\subseteq \FK{U(\g)}_r$. 
Recall that we get an isomorphism from $\FK{G_E}$ to $\FK{U(\g)}$, because
$\FK{G_E}$ is isomorphic to $\FK{M}$, which is in turn isomorphic to
$\FK{U(\g)}$ by Theorem 2.14 of \cite{M4}. This isomorphism is extended by the isomorphism from $\FK{G_E}_r$ to $\FK{U(\g)}_r$ given
in part (c) of the theorem. 

(b) Suppose that every non-zero diagonalizable one-parameter element of $E$ is $(\pm)$-\-di\-ag\-onal\-izable.
Then (\ref{eaction}) extends to an action of $\g$ if and only if (\ref{eaction2}) extends to an action of $\g^{op}$. (These actions are related
by $(x_\lo f)^{inv} =-x_\ro (f^{inv})$ where $x\in \g$, $f\in\FK{G_E}_r$, and $inv$ denotes the comorphism of the inverse map of $G_E$ .) 
Similarly, (\ref{Vaction}) extends to an action of $G_E$ if and only if 
(\ref{Vaction2}) extends to an action of $G_E^{op}$. (These actions are related
by $S^*(g_\lo h) =(g^{-1})_\ro (S^*(h))$ where $g\in G_E$, $h\in \FK{U(\g)}_r$, and $S^*$
denotes the dual of the antipode $S$ of $U(\g)$.)  
\end{Remark}

\Proof To (a): Let $f\in\FK{G_E}_r$, $e_1,\,\ldots,\, e_p\in E\setminus
\{0\}$, and $p\in\N$. Immediately from Proposition \ref{eact} (b) follows
$\left((e_1\cdots e_p)_\ro  f\right)(1)=\left((e_1\cdots e_p)_\lo f\right)(1)$. Since $U(\g)$ is spanned
by $1$ and the elements $e_1 \cdots e_p$, where $e_1,\,\ldots,\, e_p\in E\setminus
\{0\}$, and $p\in\N$, this shows the second equality in (\ref{phidef}).

Let $\Phi:\FK{G_E}_r\to U(\g)^*$ be the linear map defined by (\ref{phidef}). Then $\Phi$ is $U(\g)$-equivariant because of
\begin{eqnarray*}
  \left(x_\ro\Phi(f)\right)(y)=\Phi(f)(yx)=\left((yx)_\ro f\right)(1)=  \left(y_\ro (x_\ro f)\right)(1)=\left(\Phi (x_\ro f)\right)(y) .
\end{eqnarray*}
for all $x,y\in U(\g)$, $f\in\FK{G_E}_r$. Similarly follows the $U(\g)^{op}$-equivariance of $\Phi$. 

The elements of $\g$ act as derivations on $\FK{G_E}_r$, since the elements
$e\in E\setminus\{0\}$ act as derivations. We find $\Phi(1)=1$. The linear map $\Phi$ is also multiplicative:
We have
\begin{eqnarray*}
   \left(\Phi(f)\cdot\Phi(\ti{f})\right)(1)=
   (\Phi(f)\otimes\Phi(\ti{f}))(1\otimes 1)=\left(\Phi(f)(1)\right)\left(\Phi(\ti{f})(1)\right)=f(1)\ti{f}(1)=
    \Phi(f\ti{f})(1)
\end{eqnarray*}
for all $f,\ti{f}\in\FK{G_E}_r$. Now let $a\in U(\g)$ such that 
\begin{eqnarray*}
   \left(\Phi(f)\cdot\Phi(\ti{f})\right)(a)= \Phi(f\ti{f})(a)
\end{eqnarray*}
for all $f,\ti{f}\in\FK{G_E}_r$. Then for all $x\in\g$, and $f,\ti{f}\in\FK{G_E}_r$ we find
\begin{eqnarray*}
   \left(\Phi(f)\cdot\Phi(\ti{f})\right)(ax)
   =\left(x_\ro(\Phi(f)\cdot\Phi(\ti{f})) \right)(a)
   =\left((x_\ro \Phi(f))\cdot\Phi(\ti{f}) +\Phi(f)\cdot (x_\ro \Phi(\ti{f}))\right)(a)\\
   = \left( \Phi(x_\ro f)\cdot\Phi(\ti{f}) +\Phi(f)\cdot   \Phi(x_\ro\ti{f})\right)(a)
   = \Phi((x_\ro f) \ti{f})(a) +\Phi(f (x_\ro\ti{f}))(a)
   =\Phi(x_\ro (f \ti{f}))(a)\\
   =\left(x_\ro\Phi(f \ti{f})\right)(a)=\Phi(f \ti{f})(ax).
\end{eqnarray*}

To show that $\Phi$ is injective, and that its image is contained in $\FK{U(\g)}_r$, it is sufficient to show: Let $f\in\FK{G_E}_r$. 
Let $e_1,\,\ldots,\,e_p\in E\setminus\{0\}$, $p\in\N$. Develop
\begin{eqnarray}\label{Taylorallg1}
  f(\kappa_{e_1}(t_1)\cdots\kappa_{e_p}(t_p))=\sum_{n_1,\,\ldots,\,n_p}c_{n_1\cdots n_p}t_1^{n_1}\cdots
  t_p^{n_p}\quad\mb{ for all }\quad t_i\in\F_{e_i},\;i=1,\ldots,\,p.
\end{eqnarray}
Then we have
\begin{eqnarray}\label{Taylorallg2}
  \Phi(f)\circ  \rho_{(e_1,\,\ldots,\,e_p)}=\sum_{n_1,\,\ldots,\,n_p}c_{n_1\cdots n_p}\eta_1^{n_1}\otimes\cdots \otimes \eta_p^{n_p}.
\end{eqnarray}

To prove this equation it is sufficient to show
\begin{eqnarray}\label{Taylorallg3}
   \left({e_1^{k_1}}_\ro\left(\cdots \left({e_p^{k_p}}_\ro  f\right)\cdots\right)\right)(1) 
   =\sum_{n_1,\,\ldots,\,n_p}c_{n_1\cdots n_p}\eta_1^{n_1}(e_1^{k_1}) \cdots \eta_p^{n_p}(e_p^{k_p})
\end{eqnarray}
for all $k_1,\,\ldots,\,k_p\in\Nn$. But this follows from part (b) of Theorem \ref{eact}.\vspace*{1ex}

To (b): Let $h\in \FK{U(\g)}_r$. Let $e_1,\,\ldots,\,e_p\in E\setminus\{0\}$, $p\in\N$, and $g_i\in V_{e_i}$, $i=1,\,\ldots,\,p$. 
Immediately from Theorem \ref{Vact} (b) follows $\left((g_1\cdots g_p)_\ro h\right)(1)=\left((g_1\cdots g_p)_\lo h\right)(1)$. This
shows the second equality in (\ref{xidef}).\vspace*{1ex}

For every $h\in \FK{U(\g)}_r$ let $\Xi(h)$ be the function on $G_E$ defined by (\ref{xidef}). To show that $\Xi$ is injective, and that its 
image is contained in $\FK{G_E}_r$, it is sufficient to show: Let $h\in\FK{U(\g)}_r$. Let $e_1,\,\ldots,\,e_p\in E\setminus\{0\}$, $p\in\Nn$. 
Develop
\begin{eqnarray}\label{uTaylor1}
    h\circ\rho_{(e_1,\,\ldots,\, e_p)}=\sum_{n_1,\,\ldots,\,n_p}c_{n_1\cdots n_p}\eta_1^{n_1}\otimes\cdots\otimes \eta_p^{n_p}.
\end{eqnarray}
Then we have
\begin{eqnarray}\label{uTaylor2}
    \Xi(h)\left(\kappa_1(t_1)\cdots\kappa_p(t_p)\right)   =\sum_{n_1,\,\ldots,\,n_p}c_{n_1\cdots n_p} t_1^{n_1}\cdots t_p^{n_p}\quad
    \mb{ for all }\quad t_i\in\F_{e_i},\;i=1,\ldots,\,p.
\end{eqnarray}
But this follows by the definition of $\Xi$ and part (b) of Theorem \ref{Vact}.\vspace*{1ex} 

The proof that $\Xi$ is a $G_E^{op}\times G_E$-equivariant linear map, is similar to the proof that the map $\Phi$ of part (a) is a 
$U(\g)^{op}\otimes U(\g)$-equivariant linear map. It remains to show that $\Xi$ is a morphism of algebras.
It is not difficult to check $\Xi(1)=1$ directly, or it is possible to use
that formula (\ref{uTaylor1}) implies (\ref{uTaylor2}). Now let $h,\,\ti{h}\in \FK{U(\g)}_r$. Let 
$\underline{e}=(e_1,\,\ldots,\,e_p)\in (E\setminus\{0\})^p$, $p\in\Nn$. Develop
\begin{eqnarray*}
    \rho_{\underline{e}}^* (h) =\sum_{\underline{n}} c_{\underline{n}}(h)\,\eta^{\underline{n}}
     \quad\mb{ and }\quad \rho_{\underline{e}}^*(\ti{h})=\sum_{\underline{n}}
     c_{\underline{n}}(\ti{h})\,\eta^{\underline{n}}\quad\mb{ and }\quad
   \rho_{\underline{e}}^*(h\ti{h})=\sum_{\underline{n}}
     c_{\underline{n}}(h\ti{h})\,\eta^{\underline{n}},
\end{eqnarray*}
where $\underline{n}:=(n_1,\,\ldots,\,n_p)$ and $\eta^{\underline{n}}:=\eta_1^{n_1}\otimes \cdots \otimes \eta_p^{n_p}$. 
Note also that only finitely many of the coefficients $c_{\underline{n}}(h)$, $c_{\underline{n}}(\ti{h})$, $c_{\underline{n}}(h\ti{h})$ 
are non-zero. By Proposition \ref{ureg} (b) we have 
$\rho_{\underline{e}}^*(h\ti{h})=\rho_{\underline{e}}^*(h)\bullet\rho_{\underline{e}}^*(\ti{h})$. It follows
\begin{eqnarray*}
   c_{\underline{k}}(h\ti{h})=\sum_{\underline{m},\,\underline{n},\;\underline{m}+\underline{n}=\underline{k}}
     c_{\underline{m}}(h)c_{\underline{n}}(\ti{h})\quad\mb{ for all }\quad\underline{k}.
\end{eqnarray*}
Since (\ref{uTaylor1}) implies (\ref{uTaylor2}) from the last equations follows
\begin{eqnarray*}
   \Xi(h\ti{h})(\kappa_{\underline{e}}(\underline{t})) \,=\,
 \Xi(h)(\kappa_{\underline{e}}(\underline{t})) \; \Xi(\ti{h})(\kappa_{\underline{e}}(\underline{t})) 
 \quad \mb{ for all }\quad \underline{t}\in\F_{\underline{e}}.
\end{eqnarray*}

To (c): This follows because (\ref{Taylorallg1}) implies (\ref{Taylorallg2}), and  (\ref{uTaylor1}) implies (\ref{uTaylor2}).
\qed

As the proof shows, the embeddings of the last theorem can be described nicely
if $E$ consists only of locally nilpotent one-parameter elements. We 
state this as a Corollary.
Note also that as a particular case of Theorem \ref{eact} (b) we get the following Taylor formula for regular functions. It is valid
independently of any assumptions concerning the actions of $\g$ and $\g^{op}$ on $\FK{G_E}_r$. Only the derivations (\ref{eaction}) and
(\ref{eaction2}) are applied several times.

\begin{Corollary} Suppose that $E$ consists only of locally nilpotent one-parameter elements. Then for a regular function $f\in\FK{G_E}_r$ the
  following Taylor formula holds:
\begin{eqnarray*}
  f(\exp(t_1 e_1)\cdots \exp(t_p e_p) ) = \sum_{k_1,\,\ldots,\, k_p\in\Nn}c_{k_1\cdots k_p}\,\frac{t_1^{k_1}\cdots t_p^{k_p}}{k_1!\cdots k_p!}
\end{eqnarray*} 
with
\begin{eqnarray*}
  c_{k_1\cdots k_p}= \left({e_1^{k_1}}_\ro(\cdots ({e_p^{k_p}}_\ro f)\cdots)\right)(1)= 
                    \left({e_p^{k_p}}_\lo(\cdots ({e_1^{k_1}}_\lo f)\cdots)\right)(1),
\end{eqnarray*}
where $e_1,\,\ldots,\, e_p\in E\setminus\{0\}$, $t_1,\,\ldots,\,t_p\in\F$, and $p\in\Nn$.
\end{Corollary}
Now the embeddings are described by reading the Taylor formula from the left to the right resp. from the right to the left.
\begin{Corollary}\label{nilmor} Suppose that $E$ consists only of locally nilpotent one-parameter elements.

(a) Suppose that there is an action of the Lie algebra $\g$ on $\FK{G_E}_r$ extending the action (\ref{eaction}) 
of the elements $e\in E\setminus\{0\}$. Then the embedding $\Phi:\FK{G_E}_r\to \FK{U(\g)}_r$ of Theorem \ref{RegDual} (a) is determined by
\begin{eqnarray*}
 f(\exp(t_1 e_1)\cdots \exp(t_p e_p) ) = \sum_{k_1,\,\ldots,\,
 k_p\in\Nn}\left(\Phi(f)\right)(e_1^{k_1}\cdots e_p^{k_p})\,\frac{t_1^{k_1}\cdots t_p^{k_p}}{k_1!\cdots k_p!}
\end{eqnarray*}
where $f\in\FK{G_E}_r$, $e_1,\,\ldots,\, e_p\in E\setminus\{0\}$, $t_1,\,\ldots,\,t_p\in\F$, and $p\in\Nn$.

(b)  Suppose that there is an action of the group $G_E$ on $\FK{U(\g)}_r$ extending the action (\ref{Vaction}) 
of the groups $V_e$, $e\in E\setminus\{0\}$. Then the embedding $\Xi:\FK{U(\g)}_r\to \FK{G_E}_r$ of Theorem \ref{RegDual} (b)
is determined by
\begin{eqnarray*}
 \left(\Xi(h)\right)(\exp(t_1 e_1)\cdots \exp(t_p e_p) ) = \sum_{k_1,\,\ldots,\,
 k_p\in\Nn} h(e_1^{k_1}\cdots e_p^{k_p})\,\frac{t_1^{k_1}\cdots t_p^{k_p}}{k_1!\cdots k_p!}
\end{eqnarray*}
where $h\in \FK{U(\g)}_r$, $e_1,\,\ldots,\, e_p\in E\setminus\{0\}$, $t_1,\,\ldots,\,t_p\in\F$, and $p\in\Nn$.
\end{Corollary}
%
%
%
%

%
%
%%%%%%%%%%%%%%%%%%%%%%%%%%%%%%%%%%%%%%%%%%%%%%%%%%%%%%%%%%%%%%%%%%%%%%%%%%%%%%%%%%%%%%%%%%%%%%%%%%%%%%%%%%%%%%%%%%%%%%%%%%%%%%%%%%%%%%%%%%%%%%%
%
\section{Coordinate rings of matrix coefficients of integrable representations\label{mop}}
%
%%%%%%%%%%%%%%%%%%%%%%%%%%%%%%%%%%%%%%%%%%%%%%%%%%%%%%%%%%%%%%%%%%%%%%%%%%%%%%%%%%%%%%%%%%%%%%%%%%%%%%%%%%%%%%%%%%%%%%%%%%%%%%%%%%%%%%%%%%%%%%%%
%
%
%
In this section we investigate and describe various coordinate rings associated to categories of integrable representations of Lie algebras. 
For linear algebraic integrable Lie algebra some of these coordinate rings have been already defined in \cite{K1} on the associated 
groups, but not further investigated. We work here in a more general situation, adapted to the Tannaka reconstruction. 
By using integrable representations of free Lie algebras we show, how some of the results on the shuffle algebra and on the algebra of polynomial 
functions on the free Kac-Moody group, which have been obtained by Y. Billig and A. Pianzola in \cite{BiPi}, fit into the context of the Tannaka 
reconstruction.\vspace*{1ex}

Let $\g$ be a Lie algebra. Fix $E_d,\,E_{ln}\subseteq \g$ such that $ E:=E_d \cup E_{ln}$ generates $\g$.

\begin{Definition}\label{R1} A $\g$-module $V$ is called integrable (with respect to $E_d$ and $E_{ln}$) if the elements of
  $E_d$ act diagonalizable by integer eigenvalues on $V$, and the elements of $E_{ln}$ act locally nilpotent on $V$.
\end{Definition}

It is easy to check that submodules, quotients, sums, and tensor products of integrable $\g$-modules are again integrable.
For any $\g$-module $V$ there exists a maximal integrable $\g$-submodule $V_{int}$ of $V$, which we call the {\it integrable part} of $V$. 
(Take the sum over all integrable submodules of $V$. It is non-empty since $\{0\}$ is integrable.) If $V$, $W$ are
$\g$-modules then any $\g$-equivariant linear map $\al:V\to W$ restricts to a
$\g$-equivariant linear map $\al_{int}:V_{int}\to W_{int}$.
\subsection{\label{intfull}Integrable representations with full duals}
{\bf Assumption RR1:} In this subsection we assume that $\g$ acts faithfully on the category ${\mathcal C}_{rr}$ of integrable $\g$-modules. 
\vspace*{1ex}

Then the categories ${\mathcal C}_{rr}$ and $({\mathcal C}_{rr})^{full}$ satisfy our requirements for the Tannaka reconstruction. As a particular case of 
Theorem \ref{gone} we obtain:
\begin{Theorem}\label{R2} The pair ${\mathcal C}_{rr}$, $({\mathcal C}_{rr})^{full}$ is very good for integrating $\g$.
Let $M_{rr}$ be the corresponding Tannaka monoid. The group $G_{rr}$ generated by
\begin{eqnarray*}
  t^e \;\mb{ where }\; e\in E_d\setminus\{0\},\; t\in \F^\times,  \quad\mb{ and }\quad  \exp(te) \;
  \mb{ where }\; e\in E_{ln}\setminus\{0\},\; t\in \F,
\end{eqnarray*}
is dense in $M_{rr}$. 
\end{Theorem}

In the following we describe the coordinate ring $\FK{M_{rr}}$, resp. its isomorphic restriction $\FK{G_{rr}}$ onto $G_{rr}$.\vspace*{1ex}
Denote by $\FK{U(\g)}_{rr}$ the algebra of matrix coefficients of ${\mathcal C}_{rr}$, $({\mathcal C}_{rr})^{full}$ on $U(\g)$.
By Theorem 2.14 of \cite{M4} we get a $U(\g)^{op}\otimes U(\g)$-equivariant isomorphism of algebras 
$\Psi:\,\FK{M_{rr}}\to \FK{U(\g)}_{rr}$ by 
\begin{eqnarray*}
     \Psi(f)(x)= (x_\ro f)(1)=(x_\lo f)(1) \quad \mb{ where }\quad  f\in\FK{M_{rr}},\;x\in U(\g),
\end{eqnarray*}
resp. by
\begin{eqnarray*}
 \Psi(f_{\phi v})=g_{\phi v}\quad \mb{ where }\quad \phi\in V^*,\;v\in V,\;V\mb{ integrable }. 
\end{eqnarray*}
\begin{Theorem}\label{R3} The algebra of matrix coefficients $\FK{U(\g)}_{rr}$ of ${\mathcal C}_{rr}$, $({\mathcal C}_{rr})^{full}$ on $U(\g)$ 
is the integrable part of the $\g$-module $U(\g)^*$.
\end{Theorem}

\Proof For this proof denote the integrable part of the $\g$-module $U(\g)^*$ by $U(\g)^*_{rr}$. The inclusion 
$\FK{U(\g)}_{rr}\subseteq U(\g)^*_{rr}$ follows immediately, since the functions of $\FK{U(\g)}_{rr}$ are matrix coefficients
$g_{\phi v}$ with $\phi\in V^*$, $v\in V$, and $V$ integrable.

To show the inclusion $U(\g)^*_{rr}\subseteq \FK{U(\g)}_{rr}$ fix an element
$h\in U(\g)^*_{rr}$. The $\g$-submodule $V:=U(\g)_\ro h$ of $U(\g)_{rr}^*$ is integrable. 
Choose a base $(x_j)_\ro h$, $x_j\in U(\g)$, $j\in J$, of $V$. Let $\phi_j\in V^*$,
$j\in J$, such that $\phi_i((x_j)_\ro h)=\delta_{ij}$, $i,j\in J$. For $x\in U(\g)$
we can develop 
\begin{eqnarray*}
   x_\ro h=\sum_j c_j(x)\, (x_j)_\ro h
\end{eqnarray*} 
with coefficients $c_j(x)\in\F$, $c_j(x)\neq 0$ for only finitely many $j\in J$. Set
\begin{eqnarray*}
    \phi:=\sum_j h(x_j)\phi_j\in V^*,
\end{eqnarray*}
where the possibly infinite sum is interpreted in the obvious way.
 For all $x\in U(\g)$ we find
\begin{eqnarray*}
  g_{\phi h}(x)=\phi(x_\ro h)=\sum_j c_j(x) h(x_j)=\left(\sum_j c_j(x)\, (x_j)_\ro h\right)(1)= (x_\ro h)(1)=h(x).
\end{eqnarray*}
\qed

Recall the coordinate ring $\FK{G_{rr}}_r$ of regular functions on $G_{rr}$ and its properties from Definition \ref{DefR} and Proposition 
\ref{PropR}. 

\begin{Definition}\label{R6}
A function $f\in\FK{G_{rr}}_r$ is called restricted regular if
\begin{eqnarray*}
   \kappa_{\underline{e}}^*\left((G_{rr}^{op})_\lo f\right)
\end{eqnarray*} 
spans a finite-dimensional linear subspace of $\FK{\F_{\underline{e}}}$ for all $\underline{e}\in (E\setminus\{0\})^p$ and $p\in\N$.
Denote by $\FK{G_{rr}}_{rr}$ the set of restricted regular functions.
\end{Definition}

For an integrable $\g$-module $V$, $v\in V$, and $\phi\in V^*$ we denote the
restriction of the matrix coefficient $f_{\phi v} :\,M_{rr} \to \F$ onto
$G_{rr}$ also by $f_{\phi v}$.
\begin{Proposition}\label{R7} 
(a) $\FK{G_{rr}}_{rr}$ is a coordinate ring without zero divisors. Left and right multiplications with elements of $G_{rr}$ induce 
comorphisms of $\FK{G_{rr}}_{rr}$.

(b) Let $V$ be an integrable $\g$-module and $V^*$ the full dual. The matrix coefficients of $V$ and $V^*$ on $G_{rr}$ are contained 
in $\FK{G_{rr}}_{rr}$.
\end{Proposition}

\Proof To (b): Let $v\in V$ and $\phi\in V^*$. Let $p\in \N$, $\underline{e}\in (E\setminus\{0\})^p$. Since the elements of $E_d$ act
diagonalizable by integer eigenvalues on $V$, and the elements of $E_{ln}$ act locally nilpotent on $V$, there exist elements $v_{k_1\cdots  k_p}\in V$, 
$v_{k_1\cdots k_p}\neq 0$ for at most finitely many $k_1,\,\ldots,\,k_p$, such that 
\begin{eqnarray*}
\kappa_{\underline{e}}(t_1,\,\ldots,\,t_p)v =\sum_{k_1,\,\ldots ,\,k_p} v_{k_1\cdots k_p}\,t_1^{k_1}\cdots t_p^{k_p} \quad\mb{ for all }\quad
(t_1,\,\ldots,\,t_p)\in \F_{\underline{e}}.
\end{eqnarray*}
It follows
\begin{eqnarray*}
 (g_\lo f_{\phi v})\left(\kappa_{\underline{e}}(t_1,\,\ldots,\,t_p)\right) =\sum_{k_1,\,\ldots ,\,k_p}
\phi(g v_{k_1\cdots k_p})\,t_1^{k_1}\cdots t_p^{k_p} \quad\mb{ for all }\quad g\in G_{rr},\; (t_1,\,\ldots,\,t_p)\in \F_{\underline{e}}.
\end{eqnarray*}
This shows that $f_{\phi v}:\,G_{rr}\to\F$ is restricted regular.
\vspace*{1ex}

To (a): It is easy to check that $\FK{G_{rr}}_{rr}$ is an algebra, and the left and right multiplications with elements of $G_{rr}$ induce 
comorphisms of $\FK{G_{rr}}_{rr}$. There are no zero divisors of $\FK{G_{rr}}_{rr}$,
because it is a subalgebra of $\FK{G_{rr}}_r$, which has no zero divisors by
Proposition \ref{PropR}. It remains to show that $\FK{G_{rr}}_{rr}$ is point separating on $G_{rr}$. This follows from (b),
and from the definition of $G_{rr}$, since the full dual 
$V^*$ of an integrable
$\g$-module $V$ is point separating on $V$.
\qed

We call a $G_{rr}$-module $V$ {\it differentiable} if for every $v\in V$ and for every $e\in E\setminus\{0\}$ the following holds:
\begin{itemize}
\item If $e\in E_d\setminus\{0\}$ then $t^e v$ is Laurent polynomial in $t\in\F^\times$ with coefficients in $V$. 
\item If $e\in E_{ln}\setminus\{0\}$ then $exp(te)v$ is polynomial in $t\in\F$ with coefficients in $V$. 
\end{itemize}

For a differentiable $G_{rr}$-module $V$, $v\in V$, and $\phi\in V^*$ define the matrix coefficient $f_{\phi v}:G_{rr}\to\F$ by 
$f_{\phi v}(g):=\phi(gv)$, $g\in G_{rr}$.

\begin{Theorem}\label{R8a} The algebra of restricted regular functions $\FK{G_{rr}}_{rr}$ coincides
  with the algebra of matrix coefficients of differentiable $G_{rr}$-modules and their full duals.
\end{Theorem}

\Proof (a) We first show that $G_{rr}$ acts differentiable on every $G_{rr}$-invariant subspace $V$ of $\FK{G_{rr}}_{rr}$.
Let $e\in E_{ln}\setminus\{0\}$. Let $f\in V$. Since $f$ is restricted regular, there exists a positive integer $N\in\N$, 
such that for all $g\in G_{rr}$ we have 
\begin{eqnarray*}
   (\exp(te)_\ro f)(g)= f(g\exp(te))= (g_\lo f)(\exp(te))=\sum_{j=0}^N  a_j(g) \,t^j \quad\mb{ where }\quad a_j(g)\in \F,\;t\in\F.
\end{eqnarray*}
The linear system of equations
\begin{eqnarray*}
  \exp(m e)_\ro f =\sum_{k=0}^N m^k a_k \quad \mb{ where }\quad m=0,1,\ldots,\,N,
\end{eqnarray*}
which is a Van-der-Monde system, can be solved for $a_0,\,\ldots,\, a_N$, showing that these functions are $\Q$-linear combinations of the 
functions $\exp(0 e)_\ro f,\,\ldots,\,\exp(N e)_\ro f $ contained in $V$.  In particular the functions $a_0,\,\ldots,\,a_N$ belong to $V$.
 
Similarly, if $e\in E_d\setminus\{0\}$ and $f\in V$ then $(t^e)_\ro f$ is Laurent polynomial in $t\in\F^\times$ with coefficients in $V$.
\vspace*{1ex}

(b) It is easy to check that every matrix coefficient of a differentiable $G_{rr}$-module and its full dual is contained in $\FK{G_{rr}}_{rr}$. 
Now let $f\in \FK{G_{rr}}_{rr}$. Let $V$ be the $\F$-linear subspace of $\FK{G_{rr}}_{rr}$ generated by
the functions $(G_{rr})_\ro f$. By part (a) of this proof the $G_{rr}$-module $V$ is differentiable. 
Now we proceed in a similar way as in the proof of Theorem \ref{R3}. Thin out
$(G_{rr})_\ro f$ to a base $(g_j)_\ro f$, $j\in J$, of $V$. Let $\phi_j\in V^*$, $j\in J$, 
such that $\phi_j((g_i)_\ro f)=\delta_{ij}$, $i,j\in J$.
For $g\in G_{rr}$ the function $g_\ro f$ can be written in the form $ g_\ro f=\sum_j c_j(g) \, (g_j)_\ro f$
with $c_j(g)\in\F$, $c_j(g)\neq 0$ for only finitely many $j\in J$. Set $\phi:= \sum_j f(g_j)\phi_j \in V^*$,
where the possibly infinite sum has to be interpreted in the obvious way. Then
for all $g\in G_{rr}$ we have
\begin{eqnarray*}
 f_{\phi f}(g)=\phi( g_\ro f) = \sum_j c_j(g) f(g_j)=\left(\sum_j c_j(g) \, (g_j)_\ro f\right)(1)=(g_\ro f)(1)=f(g).
\end{eqnarray*} 
\qed

If $V$ is a differentiable $G_{rr}$-module then every $e\in E\setminus\{0\}$ defines an endomorphism of $V$ by 
\begin{eqnarray}\label{intformulas}
    e v:=\left\{ \begin{array}{ccc}
                \frac{d}{dt}\res{t=1} t^e v  & \mb{ if }& e\in E_d\setminus\{0\}\vspace*{1ex}\\
                \frac{d}{dt}\res{t=0} \exp(te) v  &\mb{ if } & e\in E_{ln}\setminus\{0\}
                \end{array}\right\} \quad \mb{ where }\quad v\in V.
\end{eqnarray}
To show that $\FK{G_{rr}}_{rr}$ coincides with $\FK{G_{rr}}$ we need:\vspace*{1ex}\\  
{\bf Assumption RR2:} In the rest of this subsection we assume that every differentiable $G_{rr}$-module $V$ gets the
structure of an integrable $\g$-module by (\ref{intformulas}).

\begin{Remark}  For the groups associated to the integrable complex Lie algebras in \cite{K1}, \S1.5,
  this assumption has been formulated as a conjecture, even for analytical differentiable actions.   
\end{Remark}

If this assumption holds then the category of integrable $\g$-modules is isomorphic to the category of differentiable $G_{rr}$-modules.
By Theorem \ref{R8a} we get:

\begin{Corollary}\label{R8} The coordinate ring $\FK{G_{rr}}$ coincides with the algebra of restricted regular functions $\FK{G_{rr}}_{rr}$.
\end{Corollary}

\subsection{\label{intint}Integrable representations with point separating integrable duals}
Similarly as for a $\g$-module, we call a $\g^{op}$-module $V$ integrable if the elements of $E_d$ act
diagonalizable by integer eigenvalues on $V$, and the elements of $E_{ln}$ act
locally nilpotent on $V$.
 
If $V$ is an integrable $\g$-module we call the $\g^{op}$-integrable part  of $V^*$ the {\it integrable dual}
$V^{(*)}$ of $V$. 
Note that the integrable dual $V^{(*)}$ does not have to be point
separating on $V$. For example let $\g$ be a one-dimensional abelian Lie
algebra spanned by $E=E_{ln}=\{e\}$. Let $V$ be a linear space with base
$b_n$, $n\in\Nn$, on which $\g$ acts by
\begin{eqnarray*}
  e b_n:=b_{n-1} \quad \mb{ for } \quad n\in\N,\quad \mb{ and }\quad e b_0:=0.
\end{eqnarray*}
Define $\phi_n\in V^*$, $n\in\Nn$, by $\phi_n(b_m)=\delta_{nm}$,
$n,m\in\Nn$. The elements of $V^*$ identify in the obvious way with 
the formal sums $\sum_{n\in\Nn}c_n\phi_n$, where $c_n\in\F$. The dual map $e^*$ acts on such an
element by 
\begin{eqnarray*}
   e^*\sum_{n\in\Nn} c_n\phi_n =\sum_{n\in\Nn} c_n \phi_{n+1}.
\end{eqnarray*}
This implies that if $\phi\in V^*\setminus\{0\}$ then $(e^*)^n \phi\neq 0$
for all $n\in\Nn$. It follows $V^{(*)}=\{0\}$.

If $V$ is an integrable $\g$-module then
\begin{eqnarray*}
  V_0:=\Mklz{v\in V}{\phi(v)=0\mb{ for all } \phi\in V^{(*)} }
\end{eqnarray*}
is a $\g$-submodule. The quotient $V/V_0$ is an integrable $\g$-module. Its
integrable dual $(V/V_0)^{(*)}$ identifies with $V^{(*)}$ in the obvious way.
In particular it is point separating on $V/V_0$.\vspace*{1ex}

{\bf Assumption DRR1:} In this subsection we assume that $\g$ acts faithfully on
the category ${\mathcal C}_{drr}$ of integrable $\g$-modules $V$, for which the integrable dual
$V^{(*)}$ is point separating on $V$. \vspace*{1ex}

It is easy to check that the integrable duals of the objects of ${\mathcal C}_{drr}$ give a category of duals 
$({\mathcal C}_{drr})^{int}$ as needed for the Tannaka reconstruction. As a particular case of Theorem \ref{gone} we obtain:
\begin{Theorem}\label{R9} The pair ${\mathcal C}_{drr}$, $({\mathcal C}_{drr})^{int}$ is very good for integrating $\g$.
Let $M_{drr}$ be the associated Tannaka monoid. The group $G_{drr}$ generated by
\begin{eqnarray*}
  t^e \;\mb{ where }\; e\in E_d\setminus\{0\},\; t\in \F^\times,  \quad\mb{ and }\quad  \exp(te) \;\mb{ where }\; e\in E_{ln}\setminus\{0\},\; 
   t\in \F,
\end{eqnarray*}
is dense in $M_{drr}$. 
\end{Theorem}
In the following we describe the coordinate ring $\FK{M_{drr}}$, resp. its isomorphic restriction $\FK{G_{drr}}$ onto $G_{drr}$.\vspace*{1ex}
Denote by $\FK{U(\g)}_{drr}$ the algebra of matrix coefficients of ${\mathcal C}_{drr}$, $({\mathcal C}_{drr})^{int}$ on $U(\g)$.
By Theorem 2.14 of \cite{M4} we get a $U(\g)^{op}\otimes U(\g)$-equi\-va\-ri\-ant isomorphism of algebras $\Psi:\,\FK{M_{drr}}\to \FK{U(\g)}_{drr}$ by 
\begin{eqnarray*}
     \Psi(f)(x)= (x_\ro f)(1)=(x_\lo f)(1) \quad \mb{ where }\quad f\in\FK{M_{drr}},\;x\in U(\g),
\end{eqnarray*}
resp. by
\begin{eqnarray*}
   \Psi(f_{\phi v})=g_{\phi v} 
\end{eqnarray*}
where $V$ is an integrable $\g$-module with point separating integrable dual
$V^{(*)}$, $v\in V$, and $\phi\in V^{(*)}$.

\begin{Theorem}\label{R10} The algebra of matrix coefficients $\FK{U(\g)}_{drr}$ of ${\mathcal C}_{drr}$, $({\mathcal C}_{drr})^{int}$ on
$U(\g)$ is the $\g^{op}$ and $\g$-integrable part of $U(\g)^*$, i.e., the sum over all $\g^{op}$ and $\g$-integrable submodules of $U(\g)^*$.
\end{Theorem}

\Proof For this proof denote by $U(\g)_{drr}^*$ the $\g^{op}$ and $\g$-integrable part of $U(\g)^*$. The inclusion 
$\FK{U(\g)}_{drr}\subseteq U(\g)^*_{drr}$  follows, because the functions of $\FK{U(\g)}_{drr}$ are matrix coefficients
$g_{\phi v}$ with $\phi\in V^{(*)}$, $v\in V$, where $V$ and $V^{(*)}$ are integrable.

To show the inclusion $U(\g)^*_{drr}\subseteq\FK{U(\g)}_{drr}$ fix an element $h\in U(\g)^*_{drr}$. The $\g$-module 
$V:= U(\g)_\ro h$ is integrable. Similarly as in the proof of Theorem \ref{R3} we find an element $\phi\in V^*$ such that 
\begin{eqnarray*}
  \phi( x_\ro h)= h(x) \quad\mb{ for all }\quad x\in U(\g).
\end{eqnarray*}  
Equip $V^*$ in the usual way with the structure of a $\g^{op}$ resp. $U(\g)^{op}$-module. We find
\begin{eqnarray}\label{phixh}
 (y\phi)(x_\ro h)=\phi(y_\ro (x_\ro h))=(x_\ro h) (y) = h(yx)=(y_\lo h)(x)
\end{eqnarray}
for all $x,y\in U(\g)$. Regard the $\g^{op}$-submodule $U(\g)^{op}\phi$ of $V^*$. We get a $\g^{op}$-equivariant injective linear map
\begin{eqnarray*}
   \eta: U(\g)^{op}\phi \to U(\g)^*_{drr}
\end{eqnarray*} 
by $\eta(\psi)(x):=\psi(x_\ro h)$ where $\psi\in U(\g)^{op}\phi$ and $x\in U(\g)$. Since $U(\g)^*_{drr}$ is integrable as $\g^{op}$-module also 
$U(\g)^{op}\phi$ is integrable as $\g^{op}$-module. 

Now we show that $U(\g)^{op}\phi$ separates the points of $V$. Let
$x\in U(\g)$ and consider the element $x_\ro h$ of $V$. If 
\begin{eqnarray*}
  (y\phi)\left(x_\ro h\right)=0 \quad \mb{ for all }y\in U(\g)^{op},
\end{eqnarray*}
then by (\ref{phixh}) we get $h(yx)=0$ for all $y\in U(\g)$. This is
equivalent to $x_\ro h=0$.

Therefore we have shown $h=g_{\phi h}$ with $h\in V$, $\phi\in V^{(*)}$, $V$ integrable, and $V^{(*)}$ point separating on $V$.
\qed
Recall the coordinate ring $\FK{G_{drr}}_r$ of regular functions on $G_{drr}$ and its properties from Definition \ref{DefR} and Proposition 
\ref{PropR}.  

\begin{Definition}\label{R11}
A function $f\in\FK{G_{drr}}_r$ is called doubly restricted regular if 
\begin{eqnarray*}
   \kappa_{\underline{e}}^*\left((G_{drr}^{op})_\lo f\right)
         \quad\mb{ as well as }\quad 
   \kappa_{\underline{e}}^*\left((G_{drr})_\ro f\right)  
\end{eqnarray*}
span finite-dimensional linear subspaces of $\FK{\F_{\underline{e}}}$ for all $\underline{e}\in (E\setminus\{0\})^p$ and $p\in\N$.
Denote by $\FK{G_{drr}}_{drr}$ the set of doubly restricted regular functions.
\end{Definition}

We define differentiable $G_{drr}$-modules in the same way as for $G_{rr}$. 
For every $G_{drr}$-module $V$ there exists a maximal differentiable submodule
$V_{diff}$, the {\it differentiable part} of $V$. 
Similarly we define differentiable $G_{drr}^{op}$-modules.
If $V$ is a differentiable $G_{drr}$-module then the {\it differentiable dual} $V^{(*)}$ is the
$G_{drr}^{op}$-differentiable part of $V^*$.\vspace*{1ex}

To show that $\FK{G_{drr}}_{drr}$ coincides with $\FK{G_{drr}}$ we need the following assumption:\vspace*{1ex}\\
{\bf Assumption DRR2:} In the following part of this subsection we assume that
every differentiable $G_{drr}$-module $V$, for which the differentiable dual
separates the points of $V$, gets the
structure of an integrable $\g$-module by (\ref{intformulas}).\vspace*{1ex}

\begin{Theorem}\label{R13} The coordinate ring $\FK{G_{drr}}$ coincides with the algebra of doubly restricted regular functions 
$\FK{G_{drr}}_{drr}$. 
\end{Theorem}

\Proof (a) Let $V$ be an integrable $\g$-module such that the integrable dual $V^{(*)}$ is point separating on $V$. It follows similarly 
as in the proof of Proposition \ref{R7} (b) that the matrix coefficients of $V$ and
$V^{(*)}$ on $G_{drr}$ are contained in $\FK{G_{drr}}_{drr}$.\vspace*{1ex}

(b) In the same way as in part (a) of the proof of Theorem \ref{R8a} we find that $G_{drr}^{op}$ acts differentiable on every $G_{drr}^{op}$-invariant
subspace $V$ of $\FK{G_{drr}}_{drr}$.\vspace*{1ex}

(c) Let $f\in \FK{G_{drr}}_{drr}$. By a similar proof as for Theorem \ref{R8a} there exists an differentiable  $G_{drr}$-module $W$, there 
exist $w\in W$ and $\phi\in W^*$, such that $W$ is spanned by $G_{drr}w$, and
\begin{eqnarray*}
   \phi(gw)=f(g) \quad\mb{ for all }\quad g\in G_{drr}.
\end{eqnarray*}
Let $G_{drr}^{op}$ act on $W^*$ in the obvious way. We get
\begin{eqnarray}\label{hgweq}
   h\phi(gw)=\phi(hgw)=f(hg)=(h_\lo f)(g) \quad\mb{ for all }\quad g,h\in G_{drr}.
\end{eqnarray}
We show that $\phi$ is contained in the differentiable dual $W^{(*)}$: Let $U$ be the subspace of $W^*$ spanned by $G_{drr}^{op}\phi$. 
Because of (\ref{hgweq}) we get an injective $G_{drr}^{op}$-equivariant linear map $\eta: U\to \FK{G_{drr}}_{drr}$ by
$\eta(\psi)(g):=\psi(gw)$, where $\psi\in U$ and $g\in G_{drr}$. By part (b) of this proof $G_{drr}^{op}$ acts differentiable on the image of
$\eta$. Therfore it also acts differentiable on $U$. 

If $W^{(*)}$ is not point separating on $W$ then $W_0:=\Mklz{\ti{w}\in W}{\psi(\ti{w})=0 \mb{ for all }\psi\in W^{(*)}}$ is a $G_{drr}$-submodule 
of $W$. The space $W^{(*)}$ identifies with the point separating differentiable dual of the differentiable $G_{drr}$-module $W/W_0$ in the
obvious way. By assumption DRR2, $W/W_0$ is an integrable $\g$-module. Furthermore it follows
that its $\g^{op}$-integrable dual contains the $G_{drr}^{op}$-differentiable dual. In
particular, the $\g^{op}$-integrable dual is point separating on
$W/W_0$. Furthermore, we have
\begin{eqnarray*}
   f(g)=\phi(gw)=\phi(g(w+W_0))\quad \mb{ for all }\quad g\in G_{drr}.
\end{eqnarray*}
\qed 

\subsection{Finite-dimensional integrable representations}
It is easy to check that for a finite-dimensional integrable $\g$-module $V$ the integrable dual $V^{(*)}$ coincides with the full dual $V^*$.
\vspace*{1ex}\\
{\bf Assumption FR1:} In this subsection we assume that $\g$ acts faithfully on
the category ${\mathcal C}_{fr}$ of finite-dimensional integrable $\g$-modules.\vspace*{1ex}

As a particular case of the Theorem \ref{gone} we obtain:
\begin{Theorem}\label{R14} The pair ${\mathcal C}_{fr}$, $({\mathcal C}_{fr})^{int}=({\mathcal C}_{fr})^{full}$ is very good for integrating $\g$.
Let $M_{fr}$ be the associated Tannaka monoid. The group $G_{fr}$ generated by
\begin{eqnarray*}
  t^e \;\mb{ where }\; e\in E_d\setminus\{0\},\; t\in \F^\times,  \quad\mb{ and }\quad  
                                                       \exp(te)  \;\mb{ where }\; e\in E_{ln}\setminus\{0\},\; t\in \F,
\end{eqnarray*}
is dense in $M_{fr}$. 
\end{Theorem}
In the following we describe the coordinate ring $\FK{M_{fr}}$, resp. its isomorphic restriction $\FK{G_{fr}}$ onto $G_{fr}$.\vspace*{1ex}

As a particular case of Theorem 2.22 and Theorem 2.23 of \cite{M4} we obtain:
\begin{Theorem} $M_{fr}$ is a group. The coordinate ring $\FK{M_{fr}}$ is in  the  natural way a Hopf algebra.
The set $M_{fr}$ identifies with $\Spm\FK{M_{fr}}$. The linear space $Lie(M_{fr})$ identifies with $Der_1(\FK{M_{fr}})$.
\end{Theorem} 

Denote by $\FK{U(\g)}_{fr}$ the algebra of matrix coefficients of ${\mathcal C}_{fr}$, $({\mathcal C}_{fr})^{int}=({\mathcal C}_{fr})^{full}$ on 
$U(\g)$. By Theorem 2.14 of \cite{M4} there is a $U(\g)^{op}\otimes
U(\g)$-equi\-va\-ri\-ant isomorphism of algebras $\Psi:\,\FK{M_{fr}}\to \FK{U(\g)}_{fr}$ given by 
\begin{eqnarray*}
     \Psi(f)(x)= (x_\ro f)(1)=(x_\lo f)(1) \quad \mb{ where }\quad f\in\FK{M_{fr}},\;x\in U(\g),
\end{eqnarray*}
resp. by
\begin{eqnarray*}
 \Psi(f_{\phi v})=g_{\phi v}\quad \mb{ where }\quad \phi\in V^*=V^{(*)},\;v\in
 V,\;V\mb{ finite-dimensional and integrable }. 
\end{eqnarray*}
\begin{Theorem}\label{R15} The algebra of matrix coefficients  $\FK{U(\g)}_{fr}$ of 
${\mathcal C}_{fr}$, $({\mathcal C}_{fr})^{int}=({\mathcal C}_{fr})^{full}$ on $U(\g)$
consists of the elements $f\in  U(\g)^*$, such that $U(\g)_\ro f$ is a finite-dimensional integrable $\g$-module.
\end{Theorem}                           

\Proof For this proof denote by $U(\g)_{fr}^*$ be the set of elements $f\in U(\g)^*$, such that $U(\g)_\ro f$ is a finite-dimensional
integrable $\g$-module. The inclusion $\FK{U(\g)}_{fr}\subseteq U(\g)_{fr}^*$ follows immediately, since the functions 
of $\FK{U(\g)}$ are matrix coefficients $g_{\phi v}$ with $\phi\in V^*$, $v\in V$, and $V$ a finite-dimensional and integrable $\g$-module. 

To show the reverse inclusion fix an element $h\in U(\g)_{fr}^*$. The submodule $U(\g)_\ro h$ of $U(\g)_{fr}^*$ is integrable and finite
dimensional. Similarly as in the proof of Theorem \ref{R3} we find an element $\phi\in \left(U(\g)_\ro h\right)^*$ such that $g_{\phi h}=h$.
\qed

Recall the coordinate ring $\FK{G_{fr}}_r$ of regular functions on $G_{fr}$ and its properties from Definition \ref{DefR} and Proposition \ref{PropR}. 
\begin{Definition}\label{R16}
A regular function $f:G_{fr}\to\F$ is called finite regular if 
\begin{eqnarray*}
  (G_{fr}^{op})_\lo f \quad\mb{ as well as  }\quad (G_{fr})_\ro f
\end{eqnarray*} 
span finite-dimensional linear subspaces of $\FK{G_{fr}}_r$.
Denote by $\FK{G_{fr}}_{fr}$ the set of finite regular functions.
\end{Definition}

We define differentiable $G_{fr}$-modules in the same way as for $G_{rr}$. To show that $\FK{G_{fr}}_{fr}$ coincides with $\FK{G_{fr}}$ we need the
following assumption:\vspace*{1ex}\\
{\bf Assumption FR2:} In the following part of this subsection we assume that every finite-dimensional differentiable $G_{fr}$-module $V$ gets the 
structure of an integrable $\g$-module by (\ref{intformulas}).\vspace*{1ex}

\begin{Theorem}\label{R18} The coordinate ring $\FK{G_{fr}}$ coincides with the algebra of finite regular functions $\FK{G_{fr}}_{fr}$. 
\end{Theorem}

\Proof Let $V$ be a finite-dimensional integrable $\g$-module. Obviously the matrix coefficients of $V$ and
$V^*=V^{(*)}$ on $G_{fr}$ are contained in $\FK{G_{fr}}_{fr}$.

Let $f\in\FK{G_{fr}}_{fr}$. Similarly as in the proof of Theorem \ref{R8a} follows: The fi\-nite-\-di\-men\-sion\-al $\F$-linear subspace
generated by $(G_{fr})_\ro f$ is a differentiable $G_{fr}$-module. By Assumption FR2 it gets also the structure of an integrable
$\g$-module. Similarly as in the proof of Theorem \ref{R8a} we find an element $\phi\in ((G_{fr})_\ro f)^*$ such that $f=f_{\phi f}$.  
\qed

For a integrable $\g$-module $V$ define the support of $V$ by
\begin{eqnarray*}
     supp(V):=\Mklz{e\in E}{e_V\neq 0 }.
\end{eqnarray*}
Sometimes the category $\mathcal C_{ffr}$ of finite-dimensional integrable $\g$-modules with finite support is more important than the 
category $\mathcal C_{fr}$ of finite-dimensional integrable $\g$-modules. There
are slightly modified definitions and similar theorems as for the category
$\mathcal C_{fr}$. Because it is easy to derive, we omit it. 

\subsection{Example: Free Lie algebras} 
Let $\g\neq \{0\}$ be a free Lie algebra in the elements of a set $E$. Set $E_{ln}:=E$ and $E_d:=\emptyset$.\vspace*{1ex}

Let $\We$ be the monoid of words in the elements of $E$. We denote the empty word, which is the unit of $\We$, by $1$.
For a word $w=e_1 e_2\cdots e_k$, where $e_1,\,e_2,\,\ldots,\, e_k\in  E$, we call $l(w):=k\in\Nn$ the length of $w$, and 
$supp(w):=\{e_1,\,e_2,\,\ldots,\, e_k\}\subseteq E$ the support of $w$. We set $l(1):=0$ and $supp(1)=\emptyset$.

Recall from Chapter 0 and Chapter 1 of \cite{R} that the universal enveloping algebra $U(\g)$ identifies 
with the monoid algebra of the monoid of words $\We$ over the field
$\F$.\vspace*{1ex}

{\bf The assumptions RR1, DRR1, and FFR1:}
To show that the assumptions RR1, DRR1, and FFR1 are satisfied it is sufficient to show that $U(\g)$ acts faithfully on the finite-dimensional 
integrable $\g$-modules with finite support.

Let $x= c_1 w_1+\cdots+ c_m w_m\in U(\g)\setminus\{0\}$, where $c_i\in\F\setminus\{0\}$, and $w_i\in\We$. Set $N:=\max_{i=1}^m l(w_i)\in\Nn$, 
and $J:=\bigcup_{i=1}^m supp(w_i)$. The set $J$ is a finite subset of $E$. Let $V_N(J)$ be a $\F$-linear space with base $b_w$, 
where $w\in \We$ with $l(w)\leq N$ and $supp(w)\subseteq J$. It is finite
dimensional. The universal enveloping algebra $U(\g)$ acts on $V_N(J)$ by
\begin{eqnarray*}
  w b_{\ti{w}}:=\left\{\begin{array}{ccc}
    b_{w\ti{w}} & \mb{ if } & supp(w)\subseteq J \quad\mb{and}\quad l(w)+l(\ti{w})\leq N\\
     0          & \mb{ else }&
\end{array}\right.  .
\end{eqnarray*}
This action restricts to an integrable action of the Lie algebra $\g$ on
$V_N(J)$. The support of $V_N(J)$ is finite. Now we have
\begin{eqnarray*}
  x b_\emptyset =\sum_{i=1}^m c_i b_{w_i}\neq 0.
\end{eqnarray*}

{\bf The assumptions RR2, DRR, and FFR2:}
The assumptions RR2, DRR2, and FFR2 are trivially satisfied, because the Lie
algebra $\g$ is free in the elements of $E$.\vspace*{1ex}

{\bf The groups $G_{rr}$, $G_{drr}$, and $G_{ffr}$:}
Recall that for $e\in E$ we denote by $\F_e$ be the additive group $\F$. Recall
from Remark \ref{Gfree} that the groups $G_{rr}$, $G_{drr}$, and $G_{ffr}$ are homomorphic images of the group
\begin{eqnarray*}
  G_{free}:=\bigstar_{e\in E}\,\F_e,
\end{eqnarray*}
the homomorphisms $\kappa_{free}$ induced by the maps $\kappa_{free}(t_e):=\exp(t_e e)$, $t_e\in\F_e\subseteq G_{free}$, $e\in E$. 
To show that the groups $G_{rr}$, $G_{drr}$, and $G_{ffr}$ are isomorphic to $G_{free}$, it is sufficient to show that 
the action of $G_{free}$ on the finite-dimensional integrable $\g$-modules with finite support, which is induced by $\kappa_{free}$, is faithful. 

Every element of $G_{free}\setminus\{1\}$ can be written in the form
\begin{eqnarray*}
     t_{e_p} * t_{e_2} \cdots * t_{e_1}  
\end{eqnarray*}
where $e_1,\,e_2,\,\ldots,\,e_p\in E$ with $e_1\neq e_2\neq\cdots\neq e_p$,
and $t_{e_j}\in\F_{e_j}\setminus\{0\}$, $j=1,\,\ldots,\, p$, and $p\in\N$.
Let $V(e_1 e_2\cdots e_p)$ be a finite-dimensional $\F$-linear space with base $b_0,\,b_1,\,\ldots,\, b_p$. 
Define an action of the free Lie algebra $\g$ on $V(e_1 e_2\cdots e_p)$ by
defining the action of the generators $E$ on the base $b_0,\,b_1,\,\ldots,\, b_p$:
\begin{eqnarray*}
   e_1 b_0 := b_1,\;\;e_2 b_1 := b_2,\;\;\ldots,\;\;e_p b_{p-1} := b_p,
\end{eqnarray*} 
and in all other cases the elements of $E$ act by zero on the elements of the base. The $\g$-module $V(e_1 e_2\cdots e_p)$ is
integrable. It support is finite. Since $e_i^2$ acts by zero, $i=1,\,\ldots,\,p$, we get
\begin{eqnarray*}
    \kappa_{free}(t_{e_p} * t_{e_2} \cdots * t_{e_1} )b_0= \exp(t_{e_p} e_p)\cdots \exp(t_{e_2} e_2)\exp(t_{e_1} e_1)b_0 \;\qquad\\
    = (1+t_{e_p}e_p)\cdots (1+t_{e_2}e_2)(1+t_{e_1}e_1)b_0 \;\qquad\\
    =b_0 \;+ \;\mb{a linear combination in }b_1,\,\ldots,\, b_{p-1} \;+\;
    \underbrace{t_{e_1}t_{e_2}\cdots t_{e_p}}_{\neq 0} b_p\neq b_0
\end{eqnarray*}

The group $G_{free}$ coincides with the free Kac-Moody group ${\mathcal F}(\F)$ defined in \cite{BiPi}.\vspace*{1ex}

{\bf The coordinate rings of regular, restricted regular, and doubly restricted regular functions:}
For the rest of this subsection we identify $G_{free}$ with $G_{rr}$, $G_{drr}$, and $G_{ffr}$. The algebras of regular functions on $G_{rr}$,
$G_{drr}$, and $G_{ffr}$ coincide. 
Therefore, we have the coordinate rings
\begin{eqnarray*}
 \FK{G_{free}}_{ffr}\subseteq \FK{G_{free}}_{drr}\subseteq \FK{G_{free}}_{rr}\subseteq \FK{G_{free}}_r.
\end{eqnarray*}
Since $\g$ is free in the elements of $E$, and $G_{free}$ is the free product of the groups $V_e$, $e\in E$, also the assumptions of 
Theorem \ref{RegDual} (a) and (b), resp. Corollary \ref{nilmor} (a) and (b) are satisfied. 
Therefore these coordinate rings are $U(\g)^{op}\otimes U(\g)$-equivariant isomorphic to the algebras
\begin{eqnarray*}
   \FK{U(\g)}_{ffr}\subseteq \FK{U(\g)}_{drr}\subseteq \FK{U(\g)}_{rr}\subseteq \FK{U(\g)}_r.
\end{eqnarray*}

(a) By Proposition \ref{nilregalg} the algebra $\FK{U(\g)}_r$ consists of the set of elements $h\in U(\g)^*$ such that for all $p\in\N$, and for all 
$e_1,\,\ldots,\,e_p\in E$ we have
\begin{eqnarray*}
  h(e_1^{k_1}\cdots e_p^{k_p})\neq 0 \quad\mb{ for at most finitely many }k_1,\,\ldots,\,k_p\in\Nn.
\end{eqnarray*}
Note also that $h(e_1^{k_1}\cdots e_p^{k_p})= \left((e_1^{k_1}\cdots e_p^{k_p})_\ro h\right)(1)=\left((e_1^{k_1}\cdots e_p^{k_p})_\lo h\right)(1)$.

(b) The algebra $\FK{U(\g)}_{rr}$ consists of the set of elements $h\in U(\g)^*$ such that for all $p\in\N$, and for all $e_1,\,\ldots,\,e_p\in E$ we have
\begin{eqnarray*}
  (e_1^{k_1}\cdots e_p^{k_p})_\ro h\neq 0 \quad\mb{ for at most finitely many }k_1,\,\ldots,\,k_p\in\Nn.
\end{eqnarray*}

(c) The algebra $\FK{U(\g)}_{drr}$ consists of the set of elements $h\in U(\g)^*$ such that for all $p\in\N$, and for all $e_1,\,\ldots,\,e_p\in E$ 
we have
\begin{eqnarray*}
  (e_1^{k_1}\cdots e_p^{k_p})_\ro h\neq 0 \quad\mb{ for at most finitely many }k_1,\,\ldots,\,k_p\in\Nn,\qquad\qquad\qquad\\
  \mb{and }\quad(e_1^{k_1}\cdots e_p^{k_p})_\lo h\neq 0  \quad\mb{ for at most finitely many }k_1,\,\ldots,\,k_p\in\Nn.
\end{eqnarray*}

By Theorem \ref{RegDual}, resp. Corollary \ref{nilmor} we get an $U(\g)^{op}\otimes
U(\g)$-equivariant isomorphism of algebras $\Phi:\FK{G_{free}}_r\to \FK{U(\g)}_r$ by
\begin{eqnarray*}
  \Phi(f)(x):=(x_\ro f)(1)=(x_\lo f)(1) \quad\mb{ where }\quad f\in\FK{G_{free}}_r,\;x\in U(\g).
\end{eqnarray*}  
It restricts to the isomorphisms between $\FK{G_{free}}_{rr}$ and $\FK{U(\g)}_{rr}$, $\FK{G_{free}}_{drr}$ and $\FK{U(\g)}_{drr}$, as well as
$\FK{G_{free}}_{frr}$ and $\FK{U(\g)}_{frr}$.
Its inverse $\Phi^{-1}: \FK{U(\g)}_r\to \FK{G_{free}}_r$ is obtained by the Taylor formula
\begin{eqnarray}\label{Taylorinv}
   \Phi^{-1}(h)(\exp(t_1 e_1)\cdots \exp(t_p e_p) ) = \sum_{k_1,\,\ldots,\, k_p\in\Nn} h(e_1^{k_1}\cdots e_p^{k_p})\,
  \frac{t_1^{k_1}\cdots t_p^{k_p}}{k_1!\cdots k_p!}\,,
\end{eqnarray} 
where $h\in \FK{U(\g)}_r$, $e_1,\,\ldots,\,e_p\in E$, and $t_1,\,\ldots,\,t_p\in\F$.\vspace*{1ex}

The algebra $\FK{G_{free}}_r$ coincides with the algebra of polynomial functions $\mb{Pol}\,{\mathcal F}(\F)$ of \cite{BiPi}. The fact that the 
algebras $\FK{G_{free}}_r$ and $\FK{U(\g)}_r$ are isomorphic is about equivalent to the description of $\mb{Pol}\,{\mathcal F}(\F)$ obtained in
Theorem 2.12 (i) and Theorem 3.7 (iii) of \cite{BiPi}.\vspace*{1ex}

{\bf The algebra $\FK{U(\g)}_{ffr}$ contains, but is in general different from the shuffle algebra:}
The algebra $U(\g)^*$ can be described as follows, compare \cite{R} or check directely. 
For $w\in\We$ define 
\begin{eqnarray*}
 \phi_w\in U(\g)^*\quad  \mb{ by } \quad \phi_w(\ti{w}):=\delta_{w\,\ti{w}},\;\ti{w}\in \We.
\end{eqnarray*}
The linear space $U(\g)^*$ identifies with $\prod_{w\in{\mathcal W}} \F
\phi_w$. The algebra structure of $U(\g)^*$ is obtained by extending in the
obvious way the products 
\begin{eqnarray*}
  \phi_{w_1}\phi_{w_2}:=\sum_{(I_1,I_2)}\phi_{w(I_1,I_2)},
\end{eqnarray*}
where $w_1,\,w_2\in\We$. Here the sum runs over the tuples of sets $(I_1,I_2)$ which satisfy
$I_1\dot{\cup}I_2=\{1,\,2,\,\ldots,\,l(w)+l(\ti{w})\}$ and $|I_1|=l(w_1)$,
$|I_2|=l(w_2)$. The word $w(I_1,I_2)$ is defined by the property that the subword
determined by $I_1$ is $w_1$, and the subword determined by $I_2$ is $w_2$. It
is called a {\it shuffle} of $w_1$ and $w_2$.

The action $\pi$ of $U(\g)^{op}\otimes U(\g)$ on $U(\g)^*$ is obtained by extending
\begin{eqnarray*}
\pi(w_1\otimes w_2) \phi_w = \left\{\begin{array}{cc}
                                       \phi_{\tau} & \mb{ if there exists }
                                       \tau\in\We \mb{ such that } w=w_1\tau w_2\\
                                         0    & \mb{ else } 
                                      \end{array}\right. 
\end{eqnarray*}
where $w_1,\,w_2,\,w\in\We$.

The {\it shuffle algebra} is a $U(\g)^{op}\otimes U(\g)$-invariant subalgebra of $U(\g)^*$ defined by
\begin{eqnarray*}
   \FK{U(\g)}_{sh}:=\bigoplus_{w\in{\mathcal W}} \F\phi_w.
\end{eqnarray*}

The shuffle algebra $\FK{U(\g)}_{sh}$ is contained in $\FK{U(\g)}_{ffr}$. To show
this it is sufficient to show $\phi_w\in \FK{U(\g)}_{ffr}$ for all $w\in \We$. For
a word $w=e_1e_2\cdots e_p$, where $e_1,\,e_2,\,\ldots,\,e_p\in E$, set
\begin{eqnarray*}
  \mb{r-cut}(w)=:\{1,\;e_1,\;e_1e_2,\;e_1 e_2 e_3,\;\ldots,\; (e_1 e_2\cdots e_p)\}.
\end{eqnarray*}
Set $\mb{r-cut}(1):=\{1\}$. Then
\begin{eqnarray*}
  U(\g)_\ro\phi_w =\bigoplus_{\ti{w}\in r\mb{-}cut(w)}\F\phi_{\ti{w}}
\end{eqnarray*}
is a finite-dimensional integrable $\g$-module. By definition it follows $\phi_w\in \FK{U(\g)}_{ffr}$.

For $|E|\geq 2$ the shuffle algebra $\FK{U(\g)}_{sh}$ is a proper subalgebra of $\FK{U(\g)}_{ffr}$. 
To show this we may restrict to $E=\{e_1,e_2\}$. Let $V$ be the integrable
$\g$-module with base $b_1$, $b_2$, such that
\begin{eqnarray*}
 e_1 b_1=0,\; e_1 b_2=b_1 \quad \mb{ and }\quad e_2 b_1=b_2,\; e_2 b_2=0.
\end{eqnarray*} 
Let $\eta\in V^*$ defined by $\eta(b_1)=\eta(b_2)=1$. Then
\begin{eqnarray*}
  \FK{U(\g)}_{fr}\ni g_{\eta b_1}=\phi_1 +\phi_{e_2}
  +\phi_{e_1 e_2}+\phi_{e_2 e_1 e_2}+\phi_{e_1 e_2 e_1 e_2}+\phi_{e_2 e_1 e_2 e_ 1 e_2}+\cdots \;\notin \FK{U(\g)}_{sh}.
\end{eqnarray*}

It is well known that the shuffle algebra can be realized as the coordinate ring of the free unipotent algebraic group, compare \cite{LMa}, \cite{Ma}, 
and also \cite{Pi2}. For some different approach by free Kac-Moody groups compare \cite{Pi}, \cite{BiPi}. 
As a supplement we show that these approaches and these results, over a field of characteristic
zero, fit nicely into the situation of the Tannaka reconstruction, and can be derived easily by the methods developed so far.  

{\bf A representation theoretic interpretation of the shuffle algebra:} 
Let ${\mathcal C}_{sh}$ the category of finite-dimensional 
$\g$-modules $V$ with the following property: For every $v\in V$ there exists an integer
$N\in\Nn$, and a finite subset $J$ of $E$, such that $w v\neq 0$ at most for $w\in\We$ with $l(w)\leq N$ and $supp(w)\subseteq J$. In particular, 
such a $\g$-module $V$ is integrable and has a finite support.
Since the $\g$-modules contained in ${\mathcal C}_{sh}$ are finite-dimensional, the only possible category of duals is 
$({\mathcal C}_{sh})^{full}$.\vspace*{1ex}

The Lie algebra $\g$ acts faithfully on ${\mathcal C}_{sh}$, since it acts
faithfully on the $\g$-modules $V_N(J)$, $N\in\Nn$, $J\subseteq E$ finite, which belong to ${\mathcal C}_{sh}$.
By Theorem \ref{gone}, and by Theorem 2.23 of \cite{M4} we obtain: The pair ${\mathcal C}_{sh}$, $({\mathcal C}_{sh})^{full}$ is very good 
for integrating $\g$. The associated Tannaka monoid $M_{sh}$ is a group. The subgroup $G_{sh}$ of $M_{sh}$ generated by
\begin{eqnarray*}
      \exp(te) \quad\mb{ where }\quad e\in E,\; t\in \F,
\end{eqnarray*}
is dense in $M_{sh}$. The coordinate ring $\FK{M_{sh}}$ is isomorphic to its
restriction $\FK{G_{sh}}$ onto $G_{sh}$.\vspace*{1ex}

Similarly as above, the group $G_{free}$ is isomorphic to $G_{sh}$. The reason for this is that the $\g$-modules 
$V(e_1 e_2\ldots e_p)$, $e_1,\,e_2,\,\ldots,\,e_p\in E$, $e_1\neq e_2\cdots\neq e_p$, $p\in\N$, also belong to ${\mathcal C}_{sh}$.
Identifying $G_{free}$ with $G_{sh}$, we denote by $\FK{G_{free}}_{sh}$ the coordinate ring $\FK{G_{sh}}$.\vspace*{1ex}

By Theorem 2.14 of \cite{M4} we get a $U(\g)^{op}\otimes U(\g)$-equivariant embedding of algebras $\widetilde{\Psi}:\,\FK{M_{sh}}\to U(\g)^*$ by 
\begin{eqnarray*}
 \widetilde{\Psi}(f_{\eta v})(x):=\eta(x_V v)\quad \mb{ where }\quad \eta\in V^*,\;v\in
 V,\;V\mb{ an object of } {\mathcal C}_{sh},\;x\in U(\g). 
\end{eqnarray*}
The image of $\widetilde{\Psi}$ is $\FK{U(\g)}_{sh}$, which can be proved as follows: The
inclusion $\widetilde{\Psi}(\FK{M_{sh}})\subseteq \FK{U(\g)}_{sh}$ follows immediately from the
properties of the $\g$-modules belonging to ${\mathcal C}_{sh}$. The inclusion
$\FK{U(\g)}_{sh}\subseteq \widetilde{\Psi}(\FK{M_{sh}})$ holds, because the matrix coefficients of $V_N(J)$ give via $\widetilde{\Psi}$ the subspace
\begin{eqnarray*}
      \bigoplus_{w\in{\mathcal W},\; supp(w)\subseteq J,\;l(w)\leq N}\F\phi_w
\end{eqnarray*}
of $\FK{U(\g)}_{sh}$. Every element of $\FK{U(\g)}_{sh}$ is contained in such a
subspace for some $N\in \Nn$, and some subset $J\subseteq E$ finite.\vspace*{1ex}

In particular, from this description follows 
\begin{eqnarray*}
   \FK{G_{free}}_{sh}=\bigoplus_{w\in{\mathcal W}}\F f_w  \quad\mb{ and }\quad
   f_{w_1}f_{w_2}=\sum_{(I_1,I_2)} f_{w(I_1,I_2)},\quad w_1,\,w_2\in\We,
\end{eqnarray*} 
where as above $w(I_1,I_2)$ runs over the shuffles of $w_1$ and $w_2$. By the Taylor formula the function 
$f_w:=\widetilde{\Psi}^{-1}(\phi_w)\res{G_{free}}$ satisfies, and is determined by
\begin{eqnarray*}
  f_w\left(\exp(t_1 e_1)\cdots \exp(t_p
    e_p)\right)=\sum_{k_1,\,\ldots,\,k_p\in\Nn}   \phi_w( e_1^{k_1}\cdots e_p^{k_p})\,\frac{t_1^{k_1}\cdots t_p^{k_p} }{k_1!\cdots k_p!}
   =\sum_{k_1,\,\ldots,\,k_p\in\Nn\atop e_1^{k_1}\cdots e_p^{k_p}=w}\frac{t_1^{k_1}\cdots t_p^{k_p} }{k_1!\cdots k_p!}
\end{eqnarray*}
for all $e_1,\,\ldots,\,e_p\in E$, $p\in\N$, and $t_1,\,\ldots,\,t_p\in\F$. 

By Proposition 2.9 of \cite{BiPi} the functions $f_w$, $w\in\We$,
coincide with the functions $X_\F^w$, $w\in\We$, defined in \cite{BiPi}. 
In particular, the algebra of matrix coefficients $\FK{G_{free}}_{sh}$ on $G_{free}$
coincides with the subalgebra $X(\F)$ of $\mb{Pol}\,{\mathcal F}(\F)$ described in Theorem 3.7 (ii) of \cite{BiPi}.\vspace*{1ex} 

{\bf The group $M_{sh}$ and the Lie algebra $Lie(M_{sh})$:} By Theorems 2.22 and 2.23 of \cite{M4} the Tannaka monoid $M_{sh}$ is a group. 
It identifies as a set with $\Spm\FK{M_{sh}}$. Furthermore, $Lie(M_{sh})$ identifies as a linear space with $Der_1(\FK{M_{sh}})$. 
The algebra $\FK{M_{sh}}$ is isomorphic to the shuffle algebra, which is well investigated. In particular, the group $M_{sh}$ is described
in \cite{R}, Theorem 3.2. The Lie algebra $Lie(M_{sh})$ is described in \cite{Pi}, Theorem 4.7 (ii).

In case of a finite generating set $E$ we indicate how these descriptions can
be also obtained from the theorems developed above for the Tannaka reconstruction. 
The free Lie algebra $\g$ acts locally nilpotent on every module $V$ contained
in ${\mathcal C}_{sh}$. Recall that we denote by $\g^k$, $k\in\Nn$, the lower
central series of $\g$. Recall that for $k\in\Nn$ we defined an ideal $I_k$ of $\g$, which contains $\g^k$, by
\begin{eqnarray*}
   I_k=\Mklz{x\in\n}{x V_k=\{0\} \mb{ for all modules } V \mb{ of }{\mathcal C}_{sh}}, 
\end{eqnarray*}
where $ V_k:=\Mklz{v\in V}{x_0 x_1\cdots x_k v=0\mb{ for all } x_0,\,x_1,\,\ldots,\,x_k\in \g}$.
It is not difficult to check that the kernel of the representation $V_k(E)$ of $\g$, which has been defined above, is $\g^k$. 
Since $(V_k(E))_k=V_k(E)$ we find $I_k=\g^k$. By Theorem 2.28 of \cite{M4} the Lie algebra
$Lie(M_{sh})$ is the pro-nilpotent Lie algebra corresponding
to the lower central series $\g^k$ of $\g$, which in turn identifies
with the Lie algebra of Lie series, compare \cite{R}, Chapter 3, Section 1.
The group $M_{sh}$ is the corresponding pro-unipotent pro-algebraic group.\vspace*{1ex} 

{\bf Open questions:}
We leave it for further research to determine $M_{rr}$, $M_{drr}$,  $M_{ffr}$, and the corresponding Lie algebras in the nontrivial case $|E|>1$. (For
$|E|=1$ these monoids and their coordinate rings coincide with the additive group $\F$ equipped with its coordinate ring of polynomial functions.)

It may also be interesting to describe the Hopf algebra $\FK{U(\g)}_{ffr}$ combinatorially. Looking at the representations of $\g$ we used so far it
seems that the matrix coefficients of $\FK{G_{free}}_{ffr}$ and $\FK{G_{free}}_{sh}$ resp. $\FK{U(\g)}_{ffr}$ and $\FK{U(\g)}_{sh}$ can be related
to representations of certain quivers with and without certain sort of circles.

\subsection{Example: Kac-Moody algebras\label{Kacint}}
Let $\g$ be a symmetrizable Kac-Moody algebra over a field $\F$ of characteristic zero, associated to a generalized
Cartan matrix $A=(a_{ij})_{i,j\in I}$ with finite index set $I$. The construction of $\g$, see for example \cite{K2}, Section1, provides us
with a Cartan subalgebra $\h$ of $\g$, as well as elements $e_i$, $f_i$, which span the root spaces
$\g_{\al_i}$, $\g_{-\al_i}$, corresponding to the simple root $\al_i$, $i\in I$.
Choose a coweight lattice $H$ of $\h$, and weight lattice $P$ of $\h^*$ as described in Chapter 1 of \cite{M1}.\vspace*{1ex} 

Set $E_d:=H$ and $E_{ln}=\Mklz{e_i,\,f_i}{i\in I}$. A $\g$-module $V$ is
integrable with respect to $E_d$ and $E_{ln}$ if and only if:
\begin{itemize}
\item $V$ has a weight space decomposition $V=\oplus_{\la\in P(V)}V_\la$ with respect
  to $\h$, where the set of weights $P(V)$ is contained in $P$.
\item The elements $e_i$, $f_i$ act locally nilpotent on $V$ for all $i\in I$.
\end{itemize}

To use later, note that instead of $E_d:=H$ and $E_{ln}=\Mklz{e_i,\,f_i}{i\in I}$ equivalently we could have taken $E_d:=H$
and $E_{ln}=\bigcup_{\al\in\Delta_{re}}\g_\al$, where $\g_\al$ denotes the root space of the real root $\al\in\Delta_{re}$.
(Note also: If the generalized Cartan matrix $A$ is degenerate, then the definition of integrable which we use here is slightly more
restrictive than the definition of integrable given in \cite{KP1}, which is integrable with respect to $E_d=\emptyset$ and 
$E_{ln}=\Mklz{e_i,\,f_i}{i\in I}$. We use this definition here, because it
allows to integrate the full Kac-Moody algebra $\g$ in a reasonable way.)\vspace*{1ex}

Assumption RR1 is satisfied, since $\g$ acts already faithfully on the integrable irreducible highest weight modules. Assumption DRR1 is
satisfied since the integrable duals of the integrable irreducible highest weight modules, which coincide with the restricted duals, are point 
separating.
If $\g$ is infinite-dimensional then Assumption FR1 is not satisfied: If $\pi:\g\to End(V)$ is a finite-dimensional representation, then the
intersection of the kernel of $\pi$ with every infinite-dimensional component of $\g$ is an infinite
dimensional ideal of this component. By \cite{K2}, Proposition 1.7, which is also valid for a field of characteristic zero, 
the kernel of $\pi$ contains the derived Lie algebras of the infinite-dimensional components of $\g$.  

The (minimal) Kac-Moody group $G$, compare for its definition for example \cite{M1}, Section 1, identifies with $G_{rr}$. (The (minimal) Kac-Moody group 
defined in \cite{KP1} coincides with the derived group $G'$. Most of the theorems which hold for $G'$ can be easily adapted to
the group $G$.)
With \cite{K1}, Corollary 4.4, it follows that $G$ acts faithfully on the
integrable irreducible highest weight modules. Therefore, $G$  also identifies with $G_{drr}$. 
With \cite{K1}, Section 2, it follows that Assumptions RR2, DRR2 are satisfied.\vspace*{1ex}

The following conjecture is quite natural:
\begin{Conjecture}\label{Conjadm}
The Kac-Moody group $G$ identifies with the monoid $M_{drr}$ associated to the category
of integrable $\g$-modules, whose integrable duals are point separating, and its category of integrable duals. The Kac-Moody
algebra $\g$ is the Lie algebra of $M_{drr}$.
\end{Conjecture}

If this conjecture holds then $(G,\FK{G},\Fi{G})$ is a weak algebraic group, whose Lie algebra identifies with the Kac-Moody algebra $\g$. \vspace*{1ex}

To know $M_{drr}$ is also interesting for the following reason: It is natural to generalize the locally finite representations of a semisimple
Lie algebra by the following classes of representations of a Kac-Moody algebra: 
\begin{itemize}
\item[(1)] Sums of integrable irreducible highest weight modules.
\item[(2)] Sums of integrable irreducible lowest weight modules.
\item[(3)] Integrable modules with point separating integrable duals.
\end{itemize}
This leads to three generalizations of a semisimple simply connected linear algebraic group.
If we consider appropriate categories and their categories of integrable duals
then for possibility (1) we obtain in \cite{M1} a monoid $\GD$ whose
coordinate ring, which we now denote by $\FK{\GD}_{int}$, obeys a Peter-and-Weyl-type theorem. 
%of the form
%\begin{eqnarray*}
%  \FK{\GD}_{int}\cong\bigoplus_{\La\in P^+}L(\La)^{(*)}\otimes L(\La)
%  \quad\mb{ as }\quad \GD\times \GD\mb{-modules}.
%\end{eqnarray*}
The Lie algebra identifies by \cite{M1} with the Kac-Moody algebra $\g$. (For some more details see also the following section.)  
Up to isomorphy possibility (2) gives nothing new.   
Maybe possibility (3) is described by the conjecture.
It is remarkable that the coordinate rings of the monoids obtained by (1) and (2) resemble much more the coordinate ring of a 
semisimple simply connected linear algebraic group than the coordinate ring obtained by (3).
Perhaps these monoids are better suited to generalize some classical algebraic geometric constructions than the Kac-Moody group itself.

%It is more uncertain if the following holds:
%
%\begin{Conjecture}\label{Conjfull}
%The Kac-Moody group $G$ identifies with the monoid $M_{rr}$ associated to the category
%of integrable $\g$-modules and its category of full duals. The Kac-Moody
%algebra $\g$ identifies with the Lie algebra of $M_{rr}$.
%\end{Conjecture}

%%%%%%%%%%%%%%%%%%%%%%%%%%%%%%%%%%%%%%%%%%%%%%%%%%%%%%%%%%%%%%%%%%%%%%%%%%%%%%%%%%%%%%%%%%%%%%%%%%%%%%%%%%%%%%%%%%%%%%%%%%%%%%%%%%%%%%%%%%%%%%%
%
\section[The monoid associated to the categories ${\mathcal O}_{int}$, ${\mathcal O}_{int}^{full}$ of a Kac-Moody algebra]{The Tannaka monoid associated to the category ${\mathcal O}_{int}$ of a Kac-Moody algebra and its category of full duals\label{KMintfull}}
%
%%%%%%%%%%%%%%%%%%%%%%%%%%%%%%%%%%%%%%%%%%%%%%%%%%%%%%%%%%%%%%%%%%%%%%%%%%%%%%%%%%%%%%%%%%%%%%%%%%%%%%%%%%%%%%%%%%%%%%%%%%%%%%%%%%%%%%%%%%%%%%%
%
%
%
As in Subsection \ref{Kacint} let $\g$ be a symmetrizable Kac-Moody algebra over a field $\F$ of characteristic zero. 
The category $\mathcal O$ is defined as follows: Its objects are the $\g$-modules $V$, which have the properties:
\begin{itemize}
\item $V$ is $\h$-diagonalizable with finite-dimensional weight spaces.
\item There exist finitely many elements $\la_1,\,\ldots,\,\la_m\in\h^*$, such that the set of weights $P(V)$ of $V$ is 
contained in the union $\bigcup_{1=1}^m \Mklz{\la\in\h^*}{\la\leq \la_i}$.
\end{itemize}
The morphisms of $\mathcal O$ are the morphisms of $\g$-modules.

Define 'integrable' as in Subsection \ref{Kacint}. Let ${\mathcal O}_{int}$ be the full subcategory of the category
$\mathcal O$, whose objects are integrable $\g$-modules. This category generalizes the 
category of finite-dimensional representations of a semisimple Lie algebra, keeping the complete reducibility theorem: Every object of 
${\mathcal O}_{int}$ is isomorphic to a direct sum of the integrable
irreducible highest weight modules. 

The highest weights of the integrable irreducible highest weight modules are given by the set $P^+$ of dominant weights of the weight lattice $P$.
We fix a system of integrable irreducible highest weight modules $L(\La)$, $\La\in P^+$. Note also that $End_{\bf g}(L(\La))=\F
id_{L(\La)}$, $\La\in P^+$.

The category ${\mathcal O}_{int}$ has the properties required for the Tannaka reconstruction. In \cite{M1} we determined the monoid $\GD$ associated to the category ${\mathcal O}_{int}$ and its category of integrable duals. Its coordinate ring
$\FK{\GD}_{int}$ satisfies a Peter-and-Weyl-type theorem:
\begin{eqnarray*}
  \FK{\GD}_{int}\cong \bigoplus_{\La \in P^+}L(\La)^{(*)}\otimes L(\La) \quad \mb{ as }\quad \GD\times\GD\mb{-modules }.
\end{eqnarray*}
We showed that the Lie algebra $Lie(\GD)$ of $\GD$ can be identified with the Kac-Moody algebra $\g$. 
The Zariski-open dense unit group of $\GD$ identifies with the Kac-Mood group $G$. The restriction of the coordinate ring of $\GD$ to $G$ identifies 
with  algebra of strongly regular functions $\FK{G}_{int}$, which had already been
defined and investigated before in \cite{KP2}. Parts of the spectrum of the algebra of strongly regular functions have been investigated in
\cite{Kas}, and \cite{Pic}. The whole spectrum of $\F$-valued points of $\FK{\GD}_{int}\cong\FK{G}_{int}$ has been determined and investigated 
in \cite{M2}.\vspace*{1ex}

Instead of choosing the integrable duals it is also possible to choose the full duals. In this section we determine the monoid associated 
to the category ${\mathcal O}_{int}$ and its category of full duals. We
describe its coordinate ring of matrix coefficients. We determine its Lie algebra.
There would be the possibility to do this in a similar way as in \cite{M1}. It would be straightforward and direct, but also quite long. Instead we 
use some of the results proved in \cite{M2} for an indirect but quick approach.  

For our notation we refer to \cite{M1}, \cite{M2}, or \cite{M3}. There is only the following difference: We index the coordinate rings 
and Zariski-closures related to the category ${\mathcal O}_{int}$ and its
category of integrable duals now by '$int$'. We don't index 
the coordinate rings and Zariski-closures related to the category ${\mathcal O}_{int}$ and its category of full duals.  

\subsection{The Tannaka monoid}

\begin{Definition} We denote by $\GfD$ the Tannaka monoid associated to the category of $\g$-modules ${\mathcal O}_{int}$ and 
its category of full duals. We denote by $\FK{\GfD}$ its coordinate ring, and by $Lie(\GfD)$ its Lie algebra.
\end{Definition}

As described in Subsection 3.2 of \cite{M4} we equip $\prod_{\La\in P^+} End(L(\La))$ with the structure of a variety induced by
\begin{eqnarray}\label{EndKo1}
   \bigoplus_{\La\in P^+} \left(End(L(\La))\right)^{du}\subseteq \left(\prod_{\La\in P^+} End(L(\La))\right)^*,
\end{eqnarray}
where $\left(End(L(\La))\right)^{du}:=\left(End_{L(\La)^*}(L(\La))\right)^{du}\subseteq End(L(\La))^*$.
Specializing Theorem 3.10 and Theorem 3.7 of \cite{M4} we get:
\begin{Theorem}\label{clsub} The map 
\begin{eqnarray*} 
  \Xi_{Irr}:\;Nat \qquad\quad &\to&  \prod_{\La\in P^+} End(L(\La)) \\
   m=(m_V)_{V\;an\; obj.\;of\;{\mathcal O}_{int}} &\mapsto &\quad (m_{L(\La)})_{\La\in P^+}
\end{eqnarray*}
is an isomorphism of $\F$-algebras, which is also an isomorphism of varieties. By this map the weak algebraic monoid $\GfD$ identifies with a
closed weak algebraic submonoid of $\prod_{\La\in P^+} End(L(\La))$.
\end{Theorem}

Our first aim is to determine the monoid $\GfD$ explicitely. 
This also demonstrates some of the theorems obtained in \cite{M4}. 

\begin{Theorem}\label{RU} Fix $\al\in \rW$. For every $x\in\g_\al$ there exists an element $\exp(x)\in \GfD$, such that for every $\g$-module $V$, 
which belongs to ${\cal O}_{int}$, we have
\begin{eqnarray*}
    \exp(x)v = \exp(x_V) v \quad\mb{ for all }\quad v\in V.
\end{eqnarray*}
Furthermore, $U_\al:=\Mklz{\exp(x)}{x\in\g_\al}$ is a closed subgroup of $\GfD$. As group with coordinate ring 
it is isomorphic to the additive group $\F$ with its coordinate
ring of polynomial functions. Its Lie algebra is $Lie(U_\al)=\g_\al$.
\end{Theorem}

\Proof The elements of $\g_\al$ act locally nilpotent on the $\g$-modules contained in ${\mathcal O}_{int}$. 
By Theorem \ref{one1II} and Remark \ref{Rone1II} the theorem follows.
\qed

Set $\n_f:=\prod_{\al\in \Delta^+}\g_\al$. The Lie bracket of $\n=\oplus_{\al\in \Delta^+}\g_\al$ extends in the obvious way to a Lie
bracket of $\n_f$. For every $\g$-module $V$ contained in ${\mathcal O}_{int}$ the action of $\n$ on $V$ extends 
in the obvious way to an action of $\n_f$. Every homomorphism between $\g$-modules contained in ${\mathcal O}_{int}$ 
is also an homomorphism of $\n_f$-modules.  
Since $\n$ acts faithfully on the objects of ${\mathcal O}_{int}$ also $\n_f$
does. (To show this use the weight space decompositions of the objects of
${\mathcal O}_{int}$.) We identify $\n_f$ with its corresponding Lie subalgebra of $Nat$.
\begin{Theorem}\label{NU} For every $x\in\n_f$ there exists an element $\exp(x)\in
  \GfD$, such that for every $\g$-module $V$, which belongs to ${\cal O}_{int}$, we have
\begin{eqnarray*}
  \exp(x)v = \exp(x_V) v \quad\mb{ for all }\quad v\in V.
\end{eqnarray*}
Furthermore, $U_f:=\Mklz{\exp(x)}{x\in \n_f}$ is a closed subgroup of $\GfD$. As group with coordinate ring it is a
  pro-unipotent group. Its Lie algebra is $Lie(U_f)=\n_f$, which
  is the pro-nilpotent Lie algebra corresponding to the descending central
  series of $\n$.
\end{Theorem}

\Proof By \cite{K2}, \S 1.3, every root space $\g_\al$, $\al\in\Delta^+$, is the linear span of the
multi-brackets $[\cdots[e_{i_1},\,e_{i_2}],\,\ldots,\,e_{i_p}]$ such that
$\al_{i_1}+\al_{i_2}+\cdots+\al_{i_p}=\al$. It follows that the descending
central series of $\n$ is given by
\begin{eqnarray*}
  \n^k=\bigoplus_{\al\in\Delta^+,\;ht(\al)\geq k+1}\g_\al\quad\mb{ where }\quad k\in\Nn.
\end{eqnarray*}
It is easy to check that $\n_f$ together with the morphisms $\psi_k:\n_f\to \n/\n^k$ defined by
\begin{eqnarray*}
    \psi_k(\sum_{\al\in\Delta^+}x_\al):=(\sum_{\al\in\Delta^+,\;ht(\al)\leq k}x_\al) +\n^k,   
\end{eqnarray*}
$k\in\Nn$, is the pro-nilpotent Lie algebra corresponding to the descending central series $\n^k$, $k\in\Nn$, of the Lie algebra $\n$. 

Fix $k\in\Nn$. Recall that the linear subspace $V_k$ of an object $V$ of ${\mathcal  O}_{int}$ is defined by
\begin{eqnarray*}
       V_k:=\Mklz{v\in V}{x_0 \cdots x_k v=0 \;\mb{ for all }\; x_0,\,\ldots,\,x_k\in\n}.
\end{eqnarray*}
Recall that the ideal $I_k$ of $\n$ is defined by 
\begin{eqnarray*}
   I_k:=\Mklz{x\in\n}{x V_k=\{0\} \mb{ for all objects } V \mb{ of } {\mathcal  O}_{int}}.
\end{eqnarray*}
Furthermore, $\n^k\subseteq I^k$. If we show equality, then the theorem follows from Theorem 2.28 and Remark 2.29 of \cite{M4}.

It is easy to check that $V_k$ is $\h$-invariant for every object $V$ of ${\mathcal O}_{int}$. From this follows that also $I_k$ is $\h$-invariant 
under the adjoint action. In particular, $I_k$ is spanned by elements of root spaces.
Assume that there exists $x_\al\in (I_k\cap \g_\al)\setminus\{0\}$, where $\al\in\pW$ with $ht(\al)\leq k$. Choose
$\La\in P^{++}$ and $v_\La\in L(\La)_\La\setminus\{0\}$. Choose $y_\al\in\g_{-\al}$ such that $\iB{x_\al}{y_\al}=1$. 
By using the weight space decomposition of $L(\La)$ we find immediately $y_\al v_\La\in L(\La)_k$. Furthermore,
\begin{eqnarray*}
  x_\al y_\al v_\La= [x_\al,\,y_\al]v_\La + y_\al \underbrace{x_\al v_\La}_{=0}= \iB{\La}{\al}v_\La.
\end{eqnarray*} 
Because of $\al\in\pW$, $\La\in P^{++}$, and $\iB{\al_i}{\La}=\frac{1}{2}\iB{\al_i}{\al_i}\La(h_i)>0$, we get $\iB{\La}{\al}>0$.
In particular, $x_\al y_\al v_\La\neq 0$, which is a contradiction. 
\qed

We denote the (dual) Tits cone in $\h_\R^*$ by $X$. The sets of weights of the
$\g$-modules contained in ${\mathcal O}_{int}$ satisfy
\begin{eqnarray*}
   \bigcup_{V\;an\; object\; of\; {\mathcal O}_{int}}P(V)=X\cap P.
\end{eqnarray*}
\begin{Theorem}\label{TU} (a) For every $h\in H$ and $s\in\F^\times$ there exists an element $s^h\in\GfD$ which acts on
the $\g$-modules $V$ belonging to ${\cal O}_{int}$ by
\begin{eqnarray*}
   s^h v_\la= s^{\la(h)}v_\la\quad \mb{ for all }\quad v_\la \in V_\la,\;\la\in P(V).
\end{eqnarray*}
These elements generate a torus $T$ isomorphic to $H\otimes_\Z \F^\times$ in the obvious way.

(b) For every $\al\in Hom(X\cap P,\F)$ there exists an element
  $t(\al)\in \GfD$, such that for every $\g$-module $V$ which belongs to ${\cal O}_{int}$ we have
\begin{eqnarray*}
   t(\al)v_\la=\al(\la)v_\la\quad \mb{ for all }\quad v_\la \in V_\la,\;\la\in P(V).
\end{eqnarray*}
Furthermore, $\TD:=\Mklz{t(\al)}{\al\in Hom(X\cap P,\F)}$ is a closed submonoid of $\GfD$. As monoid with coordinate ring it identifies with
the generalized toric monoid $\widetilde{X\cap P}:=Hom(X\cap P,\F)$, $\FK{\widetilde{X\cap P}}=\FK{X\cap P}$. In particular, 
its principal open dense unit group is the torus $T$. Its idempotents are given by the elements $e(R)$, $R$ a face of $X$, which act on
the $\g$-modules $V$ belonging to ${\cal O}_{int}$ by
\begin{eqnarray*}
   e(R) v_\la=\left\{\begin{array}{ccl}
     v_\la & \mb{if} &\la\in R\\
     0     & \mb{if}  &\la\in X\setminus R
\end{array}\right.\quad \mb{ for all }\quad v_\la \in V_\la,\;\la\in P(V).
\end{eqnarray*}
The Lie algebra of $\TD$ is $Lie(\TD)=\h$.
\end{Theorem}

\Proof Follows by Theorem \ref{one2II}, and Remarks \ref{Rone2II}, and by Theorem 2.30, Remark 2.31, and Theorem 2.32 of \cite{M4}.
\qed

The following theorem is one of the main theorems of this section.

\begin{Theorem} The monoid $\widehat{G_f}$ is generated by the elements of $U_f$, by the elements of $U_\al$, $\al\in\nrW$, 
and by the elements of $\TD$. As a monoid (without coordinate ring, without Zariski topology) it
identifies by the map of Theorem \ref{clsub} with the monoid $\GfD$ introduced in \cite{M2}.
\end{Theorem}

\Proof For this proof denote by $M$ the Tannaka monoid associated to the category of $\g$-modules ${\mathcal O}_{int}$ and its category of 
full duals.  Denote by $\GfD$ be the submonoid of $M$ generated by the elements of $U_f$, of $U_\al$, $\al\in\nrW$, 
and of $\TD$.  We show $\GfD=M$.

By Theorem \ref{clsub} we identify $M$ with a closed submonoid of $\prod_{\La\in P^+}End(L(\La))$. 
We denote the closure of $A\subseteq \prod_{\La\in P^+}End(L(\La))$ as usual by
$\overline{A}$. 

By Theorem \ref{gone} the group $G$ generated by the root groups $U_\al$,
$\al\in\rW$, and the torus $T$ is dense in $M$, i.e., $\overline{G}=M$. Since the Kac-Moody group acts faithfully on the
integrable irreducible highest weight modules $L(\La)$, $\La\in P^+$, the group $G$ identifies with the Kac-Moody group. 

Now equip $\prod_{\La\in P^+}End(L(\La))$ with a second coordinate ring, generated by
\begin{eqnarray}\label{EndKo2}
   \bigoplus_{\La\in P^+} \left(End_{L(\La)^{(*)}}(L(\La))\right)^{du}\subseteq \left(\prod_{\La\in P^+} End(L(\La))\right)^*.
\end{eqnarray}
Denote the closure of $A\subseteq \prod_{\La\in P^+}End(L(\La))$ by $\overline{\overline{A}}$. Because the coordinate ring generated by 
(\ref{EndKo2}) is a subalgebra of the coordinate ring generated by (\ref{EndKo1}) it follows $\overline{A}\subseteq\overline{\overline{A}}$.

In Theorem 16 (a) of \cite{M2} we showed $\overline{\overline{G}}=\GfD$. (Please note, the closure denoted by $\overline{A}$ in \cite{M2} 
is different from the closure denoted by $\overline{A}$ here. The closure denoted by $\overline{\overline{A}}$ in \cite{M2} is the same as here.) 
It follows $\GfD\subseteq M=\overline{G}\subseteq\overline{\overline{G}}=\GfD$.
\qed

\begin{Corollary} (a) Let $G_f$ be the group generated by the elements of $U_f$, by the elements of $U_\al$, $\al\in\nrW$, and by the elements of
  $T$. Then $G_f$ is the dense open unit group of $\GfD$. It identifies with the formal Kac-Moody group.

(b) Let $G$ be the group generated by the elements of $U_\al$, $\al\in\rW$, and by the elements of
  $T$. Then  $G$ is dense in $\GfD$. It identifies with the (minimal) Kac-Moody group.
\end{Corollary}

\Proof In the proof of the last theorem we have already shown part (b). Since the formal Kac-Moody group also acts faithfully on the
integrable irreducible highest weight modules, $G_f$ identifies with the formal Kac-Moody group. Since it contains $G$, it is dense in $\GfD$. 

Obviously $G_f$ is contained in the unit group $(\GfD)^\times$. The monoid $\GfD$ is generated by $G_f$ and by the elements of $e(R)$, 
$R$ a face of the Tits cone $X$. The elements $e(R)$, $R$ different from $X$, do not act
by bijective endomorphims on every $\g$-module contained in ${\mathcal O}_{int}$. It follows $G_f=(\GfD)^\times$.

Let $\La\in P^{++}$. Let $v_\La\in L(\La)_\La\setminus\{0\}$ and $\phi_\La\in (L(\La)^{(*)})_\La\setminus\{0\}$. 
Let $\gt_\La:=f_{\phi_\La v_\La}\in\FK{\GfD}$ be the corresponding matrix
coefficient. $\La$ is contained in the interior of the Tits cone
$X$. Since the Tits cone $X$ is $\We$-invariant also $\We\La$ is contained in the
interior of $X$. The set of weights $P(\La)$ is contained in the convex hull
of $\We\La$, which lies also in the interior of $X$. Therefore, the elements
$e(R)$, $R$ a face different from $X$, act as zero on $L(\La)$. It follows
\begin{eqnarray*}
  D_{\widehat{G_f}}(\gt_\La)=D_{G_f}(\gt_\La\res{G_f}).
\end{eqnarray*}
By using the Birkhoff decomposition $G_f=\dot{\bigcup}_{n\in N}U^- n U_f$ of the
formal Kac-Moody group $G_f$ we find 
\begin{eqnarray*}
 D_{G_f}(\gt_\La\res{G_f})=U^- T U_f.
\end{eqnarray*} 
From the Bruhat or Birkhoff decompositions of the minimal Kac-Moody group $G$,
and of the formal Kac Moody group $G_f$ it follows $G_f= G U_f$. By
\cite{KP2}, Corollary 3.1, which is also valid for the slightly enlarged
Kac-Moody group $G$ which we use here, we have $G=\bigcup_{w\in{\mathcal W}} w U^- T U$,
where $\We$ denotes the Weyl group. 
It follows 
\begin{eqnarray*}
  G_f=\bigcup_{w\in {\mathcal W}} w U^-T U_f,
\end{eqnarray*}
Since left multiplications with elements of $G_f$ are Zariski-homeomorphisms, it
follows that the set $w U^-T U_f$ is open. Therefore $G_f$ is open.
\qed

There are similar structural results for $\GfD$ as for the monoid $\GD$, and the spectrum of $\F$-valued points of the coordinate ring 
$\FK{\GD}_{int}$. We do not state the results here in full length, but only give some remarks:
  
By Theorem 17 (1) of \cite{M2} there are the following Bruhat and Birkhoff decompositions 
\begin{eqnarray*}
    \GfD=\dot{\bigcup_{\hat{n}\in\widehat{N}} } U^\pm \hat{n} U_f.
\end{eqnarray*} 
Similarly as in \cite{M3} it is possible to describe the Bruhat and Birkhoff cells, as well as their closure relations.

From the Birkhoff decomposition, or by using Theorem 18 (1) of \cite{M2} it follows
\begin{eqnarray*}
 \GfD=\dot{\bigcup_{\Theta\;special}}  G e(R(\Th)) G_f.
\end{eqnarray*} 
Similarly as in \cite{M2} it is possible to define big cells of $G e(R(\Th)G_f$, describe these cells, show that countably many cover 
$G e(R(\Th)) G_f$. It is possible to find stratified transversal slices to $G e(R(\Th)) G_f$.

\subsection{The coordinate ring of matrix coefficients}
We give three descriptions of the coordinate ring of matrix coefficients $\FK{\GfD}$. They are of a similar style as the description of the coordinate
ring $\FK{G}_{int}$ obtained by V. G. Kac and D. Peterson in Theorem 1 of \cite{KP2}, and as the description of $\FK{G}_{int}$ obtained by M. Kashiwara in
Section 5 of \cite{Kas}.\vspace*{1ex}

As a particular case of Theorem 2.16 of \cite{M4} we obtain:
\begin{Theorem}\label{GfD-PW} The map 
\begin{eqnarray*}
   \bigoplus_{\La\in P^+} L(\La)^*\otimes L(\La) \;\to\;  \FK{\GfD}\quad \mb{ induced by }\quad  \phi \otimes w\mapsto f_{\phi w}
\end{eqnarray*}
is a $\GfD^{op}\times \GfD$-equivariant and $U(Lie(\GfD))^{op}\otimes U(Lie(\GfD))$-equivariant linear bijective map.
\end{Theorem}

As a particular case of Theorem \ref{R3}, the algebra of matrix coefficients $\FK{U(\g)}_{rr}$ on $U(\g)$ of the integrable $\g$-modules and
their full duals is given by
\begin{eqnarray*}
  \FK{U(\g)}_{rr} = \bigoplus_{\la\in P} \,(\FK{U(\g)}_{rr})_\la,
\end{eqnarray*}
where $(\FK{U(\g)}_{rr})_\la$ consists of the functions $l\in U(\g)^*$ which satisfy:
\begin{itemize}
\item  $h_\ro l=\la(h) l$ for all $h\in\h$.
\item  For every $i\in I$ there exists an integer $n\in\N$ such that $(e_i^n)_\ro l=(f_i^n)_\ro l=0$.
\end{itemize}

We denote the algebra of matrix coefficients of ${\mathcal O}_{int}$ and ${\mathcal O}_{int}^{full}$ on $U(\g)$ by $\FK{U(\g)}$.
By Theorem 2.14 of \cite{M4} we get a $U(\g)^{op}\otimes U(\g)$-equivariant isomorphism of algebras $\Psi:\,\FK{\GfD}\to \FK{U(\g)}$ by 
\begin{eqnarray*}
     \Psi(f)(x)= (x_\ro f)(1)=(x_\lo f)(1) \quad \mb{ where }\quad f\in\FK{\GfD},\;x\in U(\g),
\end{eqnarray*}
resp. by
\begin{eqnarray*}
 \Psi(f_{\phi v})=g_{\phi v}\quad \mb{ where }\quad \phi\in V^*,\;v\in V,\;V\mb{ an object of } {\mathcal O}_{int}. 
\end{eqnarray*}
\begin{Theorem} The algebra of matrix coefficients $\FK{U(\g)}$ is the subalgebra of $\FK{U(\g)}_{rr}$, which is given by the elements 
$l\in \FK{U(\g)}_{rr}$, such that the $\F$-linear space $U({\bf b})_\ro l$ is finite-dimensional.  
\end{Theorem}

\Proof Denote by $CR$ the subset of $\FK{U(\g)}_{rr}$ described in the theorem. With the last theorem follows, that the functions of $\FK{U(\g)}$ are
sums of matrix coefficients $g_{\phi v}$ with $\phi\in L(\La)^*$, $v\in L(\La)$, and $\La\in P^+$. 
Since $L(\La)$, $\La\in P^+$, are integrable highest weight modules, the inclusion $\FK{U(\g)}\subseteq CR$ follows.

Now let $l\in CR$. The $\g$-submodule $V:=U(\g)_\ro l$ of $\FK{U(\g)}_{rr}$ is integrable. We may write 
$V=U(\g)_\ro l=U({\bf n}^-)_\ro \left(U({\bf b})_\ro l\right)$. Here $U({\bf b})_\ro l$ is finite-dimensional. 
It follows easily that $V$ is an integrable $\g$-module of the category $\mathcal O$.
Similarly as in the proof of Theorem \ref{R3} we find an element $\phi\in V^*$, such that 
\begin{eqnarray*}
 g_{\phi l}(x)=\phi(x_\ro l)= l(x) \quad\mb{ for all }\quad x\in U(\g).
\end{eqnarray*}
\qed

Denote by $\FK{G}$ the restriction of the coordinate ring $\FK{\GfD}$ onto the Kac-Moody group $G$. Since $G$ is dense
in $\GfD$, the coordinate ring $\FK{G}$ is isomorphic to $\FK{\GfD}$ by the restriction map.
Recall the coordinate ring $\FK{G}_{rr}$, which has been introduced in Section \ref{intfull} to describe the algebra of matrix coefficients 
of the integrable $\g$-modules and its full duals.
\begin{Theorem} The coordinate ring $\FK{G}$ is the subalgebra of $\FK{G}_{rr}$ given by the functions 
$h\in \FK{G}_{rr}$, for which there exists finitely many $g_1,\,\ldots,\,g_m\in G$ such that
\begin{eqnarray*}
       u_\ro h=h  \quad \mb{ for all } \quad u\in\bigcap_{i=1}^m g_i U g_i^{-1}
 \end{eqnarray*}
\end{Theorem}
\Proof Denote by $CR$ the subset of $\FK{G}_{rr}$ described in the theorem. It is easy to check that $CR$ is a $G^{op}\times G$-invariant
subalgebra of $\FK{G}_{rr}$.

We first show $\FK{G}\subseteq CR$. By Theorem \ref{GfD-PW} it is
sufficient to show $f_{\phi v}\res{G}\in CR$ for all $\phi\in L(\La)^*$, $v\in L(\La)$, and $\La\in P^+$.
Since $L(\La)$ is integrable, we get $f_{\phi v}\res{G}\in\FK{G}_{rr}$. Since $L(\La)$ is an irreducible $G$-module, it is spanned by $G v_\La$, where
$v_\La\in L(\La)_\La\setminus\{0\}$. Therefore, $v\in L(\La)$ can be written in the form
\begin{eqnarray*}
   v=\sum_{i=1}^m c_i g_i v_\La \quad \mb{ with } \quad g_i\in G,\; c_i\in\F.
\end{eqnarray*}
For $u\in\bigcap_{i=1}^m g_i U g_i^{-1}$ it follows $uv=v$, which implies $u_\ro (f_{\phi v}\res{G})=f_{\phi v}\res{G}$.\vspace*{1ex}

Recall that the $G$-invariant subspaces, and the $\g$-invariant subspaces of $\FK{G}_{rr}$ coincide. The group $G$ acts differentiable, the Lie 
algebra $\g$ acts integrable on such a subspace, and the actions can be obtained from one another. In particular, this holds for
$CR$ and its $G$ resp. $\g$-invariant subpaces.

Next we show that $\n$ acts locally finite on $CR$. Let $h\in CR$. Then there exist elements $g_1,\,\ldots,\,g_m$ of $G$ such that $u_\ro h=h$  for all 
$u\in\bigcap_{i=1}^m g_i U g_i^{-1}$. The group $\ti{U}:=U\cap \bigcap_{i=1}^m g_i U g_i^{-1}$ is a large subgroup of $U$ in the
sense of \cite{KP2}, \S 2 B. By Lemma 2.2 of \cite{KP2} there exist finitely many real roots $\beta_1,\,\ldots,\,\beta_k$ such that
\begin{eqnarray*}
   U=U_{\beta_1}\cdots U_{\beta_k} \ti{U}.
\end{eqnarray*}
Since $G$ acts differentiable on $CR$, the $\F$-linear space $W_h$ spanned 
by $U_\ro h= \left(U_{\beta_1}\cdots U_{\beta_k}\right)_\ro h$ is finite-dimensional. 
In particular, the space $W_h$ is invariant under $U_\al=\exp(\g_\al)$ for every $\al\in \prW$. From this follows easily
that $W_h$ is invariant under $\g_\al$ for all $\al\in\prW$. Because these
root spaces generate $\n$ as a Lie algebra, $W_h$ is also invariant under $\n$.

Now let $h\in CR$. Let $V$ be the $\F$-linear subspace of $CR$ spanned by $G_\ro h$. In the proof of Theorem \ref{R8a} we have seen that there exists 
an element $\phi\in V^*$ such that 
\begin{eqnarray*}
   \phi(g_\ro h )= h(g) \quad\mb{ for all }\quad g\in G.
\end{eqnarray*} 
We may replace $V$ by $W:=U(\g)_\ro h\subseteq V$ (actually we have equality), and $\phi$ by its restriction to $W$, which we also denote by $\phi$. 
$\n$ acts locally finite on $W$. As in the proof of the last theorem we find that $W$ belongs to $\mathcal O$. In particular, 
$\GfD$ acts on $W$, extending the action of $G$. Therefore, $h\in W$ and $\phi\in W^*$ define a
matrix coefficient $f_{\phi h}\in \FK{\GfD}$ such that $f_{\phi h}\res{G}=h$.
\qed
\subsection{The Lie algebra}
Recall that $U$ is the group generated by the root grous $U_\al$, $\al\in\prW$. Denote by $\FK{U}$ the restriction of the coordinate ring
$\FK{\GfD}$ onto $U$. Denote by $\FK{U}_{int}$ the restriction of the coordinate ring $\FK{\GD}_{int}$ onto $U$.
The following theorem relates $\FK{U}_{int}$ to the coordinate ring $\FK{U_f}$
of the pro-unipotent group $U_f$. This is advantageous because
it is easy to work with pro-unipotent groups and their coordinate rings. 

\begin{Theorem} (a) The group $U$ is dense in $U_f$. In particular the coordinate
  ring $\FK{U_f}$ is isomorphic to $\FK{U}$ by the restriction map.

(b) We have $\FK{U}_{int}=\FK{U}$.
\end{Theorem}

\Proof To (a): Recall the category ${\mathcal C}(\n)$ from Section 2.8 of \cite{M4}. By the proof of Theorem 2.28 of \cite{M4} 
the group $U_f$ is the Tannaka monoid associated to ${\mathcal C}(\n)$ and ${\mathcal C}(\n)^{full}$, and $\FK{U_f}$
is its coordinate ring of matrix coefficients. The Lie algebra $\n$ is generated by $\g_\al$, $\al\in\prW$. Therefore, by Theorem \ref{gone},
the group $U$ is dense in $U_f$.

To (b): Obviously $\FK{U}_{int}\subseteq \FK{U}$. Now let $V$ be an object of ${\mathcal O}_{int}$, let $v\in V$, and $\phi\in V^*$. 
Choose a pair of dual bases
\begin{eqnarray*}
   (a_{\la i})_{\la\in P(V),\,i=1,\,\cdots,\,m_\la} \quad\mb{ and }\quad (\eta_{\la i})_{\la\in P(V),\,i=1,\,\cdots,\,m_\la}
\end{eqnarray*}
of $V$ and $V^{(*)}$ by choosing a base $a_{\la i}$, $i=1,\,\cdots,\,m_\la$, of every weight space $V_\la$, $\la\in P(V)$. Then
\begin{eqnarray*}
 \phi(u v)= \sum_{\mu,\,j}\phi(a_{\mu j})\eta_{\mu j}(uv) \quad \mb{ for all }\quad u\in U.
\end{eqnarray*}
For $w\in V$ denote by $supp(w)$ the set of weights of the non-zero weight space components
of $w$. Then $supp(uv) \subseteq (supp(v)+Q_0^+)\cap P(V)$. Since $(supp(v)+Q_0^+)\cap P(V)$ is finite, it follows 
\begin{eqnarray*}
 f_{\phi v}\res{U}= \sum_{\mu,\,j\atop\mu\in (supp(v)+Q_0^+)\cap P(V)}\phi(a_{\mu j})f_{\eta_{\mu j} v}\res{U}\,\in \FK{U}_{int}.
\end{eqnarray*}
\qed

Set $\n_f^-:=\prod_{\al\in \Delta^-}\g_\al$. Let $V$ be an object of ${\mathcal O}_{int}$, let $v\in V$ and $\phi\in V^{(*)}$. 
Every element of $\n_f^-$ maps $V=\bigoplus_{\la\in P(V)}V_\la$ into $V_f:=\prod_{\la\in P(V)}V_\la$ in the obvious way. 
Since $V^{(*)}=\bigoplus_{\la\in  P(V)}V_\la^*$ the expression
\begin{eqnarray*}
  \phi(n_-v)
\end{eqnarray*}
gives a well defined element of $\F$ for every $n_-\in\n_f^-$.

\begin{Theorem}\label{Deraosrf} We get a bijective linear map
\begin{eqnarray*}
   \n_f^-\oplus \h\oplus \n_f^+\to Der_1(\FK{\GD}_{int})
\end{eqnarray*}  
by assigning the element $x=n_-+h+n_+$, where $n_-\in\n_f^-$, $h\in\h$, and
$n_+\in\n_f^+$, the derivation $\delta_x\in Der_1(\FK{\GD}_{int})$ defined by
\begin{eqnarray*}
 \delta_x(f_{\phi v}):=\phi(n_-v)+\phi(h v)+\phi(n_+ v) \quad \mb{ for all }\quad 
  v\in V, \; \phi\in V^{(*)}, \;V \mb{ an object of } {\mathcal O}_{int}.
\end{eqnarray*}

\end{Theorem}

\Proof Choose an element $v_\La \in L(\La)_\La$, $\phi_\La \in (L(\La)^{(*)})_\La$, such that $\phi_\La(v_\La)=1$. 
Let $\gt_\La:=f_{\phi_\La v_\La}$ be the associated matrix coefficient on $\GD$. Equip the torus $T$ with its
classical coordinate ring $\FK{T}_{class}$. It identifies with the group algebra $\FK{P}=\bigoplus_{\la\in P}\F e_\la$. 
Equip the groups $U^\pm$ and the principal open set $D_{\widehat{G}}(\gt_\La)$ with the coordinate rings $\FK{U^\pm}_{int}$ and 
$\FK{D_{\widehat{G}}(\gt_\La)}_{int}$ obtained by restriction of $\FK{G}_{int}$.
By Proposition 6.4 and Theorem 6.5 of \cite{M1} the multiplication map
\begin{eqnarray*}
m: U^-\times T\times U^+\to D_{\widehat{G}}(\gt_\La)
\end{eqnarray*}
is an isomorphism of sets with coordinate rings. Denote by $j:D_{\widehat{G}}(\gt_\La)\to \GD$
the inclusion morphism.
The tangential maps of sets with coordinate rings
\begin{eqnarray*}
 \ti{T}_1(m) : Der_1(\FK{U^-}_{int})\times Der_1(\FK{T}_{class})\times
 Der_1(\FK{U^+}_{int})\to Der_1(\FK{D_{\widehat{G}}(\gt_\La)}_{int}) 
\end{eqnarray*}
and
\begin{eqnarray*}
  \ti{T}_1(j): Der_1(\FK{D_{\widehat{G}}(\gt_\La)}_{int})\to Der_1(\FK{\GD}_{int})
\end{eqnarray*}
are bijective linear maps. (Here $\ti{T}_1(m)$ is defined as concatenation with the
comorphism $m^*$, and the derivation $(\delta_1,\delta_2,\delta_3)\in Der_1(\FK{U^-}_{int})\times Der_1(\FK{T}_{class})\times
 Der_1(\FK{U^+}_{int})$ of $\FK{U^-}\otimes\FK{T}\otimes\FK{U^+}$ is given by
\begin{eqnarray*}
   (\delta_1,\delta_2,\delta_3)(f_1\otimes f_2\otimes f_3)
   =\delta_1(f_1) f_2(1)f_3(1)+f_1(1)\delta_2(f_2)f_3(1)+f_1(1)f_2(1)\delta_3(f_3)
\end{eqnarray*}
for all $f_1\in\FK{U^-}_{int}$, $f_2\in\FK{T}_{class}$, and
$f_3\in\FK{U^+}_{int}$. The map $\ti{T}_1(j)$ is defined as the concatenation
with the comorphism $j^*$.)

By the last theorem we have $\FK{U^+}_{int}=\FK{U^+}\cong \FK{U_f}$. As shown in Theorem \ref{NU} the group $U_f$ is unipotent, and its Lie algebra 
identifies  with $\n_f$. From the proof of Theorem 2.28 of \cite{M4} we know that already $Der_1(\FK{U^+})$ identifies with $\n_f$. The Chevalley 
involution $*:U^+\to U^-$ is an isomorphism of sets with coordinate ring. Therefore,
$Der_1(\FK{U^-})$ identifies with $\n_f^-$. It is well known that $Der_1(\FK{T}_{class})$ identifies with $\h$.

Let $n_-\in\n_f^-$, $h\in\h$, $n_+\in\n_f^+$ and let $\delta_{n_-}\in Der_1(\FK{U^-}_{int})$, $\delta_h\in Der_1(\FK{T}_{class})$, 
$\delta_{n_+}\in Der_1(\FK{U^-}_{int})$ be the corresponding derivations. Set
\begin{eqnarray*}
   \delta:=\ti{T}_1(j\circ m)((\delta_{n_-},\delta_h,\delta_{n_+})).
\end{eqnarray*}
Now we compute $\delta(f_{\phi v})$ for $v\in V$, $\phi\in V^{(*)}$ , $V$ an object of ${\mathcal O}_{int}$. 
Choose a pair of dual bases
\begin{eqnarray*}
   (a_{\la i})_{\la\in P(V),\,i=1,\,\cdots,\,m_\la} \quad\mb{ and }\quad (\eta_{\la i})_{\la\in P(V),\,i=1,\,\cdots,\,m_\la}
\end{eqnarray*}
of $V$ and $V^{(*)}$ by choosing a base $a_{\la i}$, $i=1,\,\cdots,\,m_\la$,
of every weight space $V_\la$, $\la\in P(V)$. Then
\begin{eqnarray*}
  f_{\phi v}(u^-tu^+)=\sum_{\la i}\phi(u^- v_{\la i})e_\la(t) \eta_{\la i}(u^+ v).
\end{eqnarray*}
for all $u^-\in U^-$, $t\in T$, and $u^+\in U^+$. Denote by $supp(v)$ the set of weights of the non-zero weight space components
of $v$. Then $\eta_{\la i}(u^+ v)$ is nonzero at most for $\la\in
(supp(v)+Q_0^+)\cap P(V)$. Since there are only finitely many such weights, it follows
\begin{eqnarray*}
   (j\circ m)^*( f_{\phi v})=\sum_{\la i\atop \la\in (supp(v)+Q_0^+)\cap P(V)}
           f_{\phi v_{\la i}}\res{U^-}\otimes e_\la\otimes  f_{\eta_{\la i} v}\res{U}.
\end{eqnarray*} 
Therefore we get
\begin{eqnarray*}
  \delta(f_{\phi v})=\sum_{\la i\atop \la\in (supp(v)+Q_0^+)\cap P(V)} 
            \left(\,\phi(n_- v_{\la i})\eta_{\la i}(v)
           +\phi(v_{\la i})\la(h)\eta_{\la i}(v)
           +\phi(v_{\la i})\eta_{\la i}(n_+ v)\,\right)\\
           =\phi(n_- v)+\phi(h v)+\phi(n_+ v).
\end{eqnarray*}
\qed

Let $\La\in P^+$. Denote by $L_{high}$ the $L(2\La)$-isotypical component of the $\g$-module
$L(\La)\otimes L(\La)$. The Kostant cone ${\mathcal V }_\La$ is defined by
\begin{eqnarray*}
  {\mathcal V }_\La:=\Mklz{v\in L(\La)}{v\otimes v\in L_{high}}.
\end{eqnarray*}
By Theorem 1 of \cite{KP1} it can be described as ${\mathcal V }_\La=G'\left( L(\La)_\La\right)$. 

\begin{Proposition}\label{Kostant} Let $\La\in P^+$. We have ${\mathcal V }_\La=\GfD  L(\La)_\La$.
\end{Proposition}

\Proof By Proposition 2.8 of \cite{M4} the $\g$-invariant subspace $L_{high}$ of $L(2\La)$ is $\GfD$-invariant. 
Therefore, the Kostant cone ${\mathcal V }_\La$ is $\GfD$-invariant. Since
$L(\La)_\La\subseteq {\mathcal V }_\La$ we get
\begin{eqnarray*}
  {\mathcal V }_\La=G'\left( L(\La)_\La\right)\subseteq \GfD L(\La)_\La\subseteq {\mathcal V }_\La.
\end{eqnarray*}
\qed

Set $\g_f:=\n^-\oplus\h\oplus\n_f$. The Lie bracket of $\g$ extends in the obvious
way to a Lie bracket of $\g_f$. Every $\g$-module contained in ${\mathcal O}_{int}$ can
be extended to a $\g_f$-module, every homomorphism between $\g$-modules
contained in ${\mathcal O}_{int}$ is also a homomorphism of $\g_f$-modules.
Since $\g$ acts faithfully on the objects of ${\mathcal O}_{int}$ also $\g_f$
does. (To show this use the weight space decompositions of the objects of ${\mathcal O}_{int}$.)
We identify $\g_f$ with its corresponding Lie subalgebra of $Nat$.

\begin{Theorem} We have $Lie(\widehat{G_f})=\g_f$.
\end{Theorem}

\Proof 
(a) We first show $\g_f\subseteq Lie(\GfD)$. By Theorem \ref{RU}, Theorem \ref{NU}, and Theorem \ref{TU} 
we have $\g_\al\subseteq Lie(\GfD)$, where $\al\in\nrW$, and $\n_f\subseteq
Lie(\GfD)$, and $\h\subseteq Lie(\GfD)$. These subalgebras of $Lie(\GfD)$ generate $\g_f$. \vspace*{1ex} 

(b) Let $\La\in P^+$. Equip $L(\La)$ with the structure of a variety induced by
$L(\La)^{(*)}\subseteq L(\La)^*$. By Lemma 3 of \cite{KP1}
the Kostant cone ${\mathcal V}_\La$ is a closed subset of $L(\La)$. We equip ${\mathcal V}_\La$ with its subvariety structure.
Fix $v_\La\in L(\La)_\La\setminus\{0\}$. Due to Proposition \ref{Kostant} we get a map
\begin{eqnarray*}
  \Phi:\;\GfD\to {\mathcal V }_\La\quad \mb{ by }\quad\Phi(m):= m v_\La, \;m\in\GfD.
\end{eqnarray*}
It is not difficult to check that $\Phi$ is a morphism of varieties.
Identify the tangent space $T_{v_\La}({\mathcal V }_\La)$ with the corresponding subspace of 
$T_{v_\La}(L(\La)_\La)$. Identify this last tangent space with $L(\La)$. 
Then the tangent map 
\begin{eqnarray*}
  T_1(\Phi ): \;Lie(\GfD) \to T_{v_\La}({\mathcal V }_\La)\quad\mb{ is given by
  }\quad \left(T_1(\Phi)\right)(x)=x_{L(\La)} v_\La,\;x\in Lie(\GfD).
\end{eqnarray*}
Furthermore, by Theorem 6.2 of \cite{M1} we have
\begin{eqnarray*}
   T_{v_\La}({\mathcal V }_\La)=\g L(\La)_\La.
\end{eqnarray*}

(c) Let $\La\in P^{++}$. Recall that the elements of $\n_f^-$ map $L(\La)=\bigoplus_{\la\in P(\La)}L(\La)_\la$ into 
$L(\La)_f:=\prod_{\la\in P(\La)}L(\La)_\la$. Let $v_\La\in L(\La)_\La$. We show that the map
\begin{eqnarray*}
   \Omega: \n_f^-\to L(\La)_f \quad\mb{ defined by }\quad \Omega(n):=n v_\La,\;n\in \n_f^- 
\end{eqnarray*}
is injective.
Let $n=\sum_{\al\in\Delta^-}n_\al\in \n_f^-$ such that $nv_\La=\sum_{\al\in\Delta^-}n_\al v_\La=0$. Then $n_\al v_\La=0$ for all $\al\in\nW$.
By Lemma 5 (b) of \cite{KP1} the set of elements of $\g$ which stabilize $L(\La)_\La$
is given by $\h\oplus\n^+$. It follows $n_\al=0$ for all $\al\in\nW$.\vspace*{1ex} 

(d) Now we show $Lie(\GfD)\subseteq \g_f$. 
Denote by $\FK{\GD}$ the restriction of $\FK{\GfD}$ onto
$\GD$. By Theorem \ref{gone} already $G$ is Zariski dense in $\GfD$. Therefore $\FK{\GfD}$ is isomorphic to $\FK{\GD}$ by the restriction map.

Let $x\in Lie(\GfD)$. Let $\ti{\delta}_x$ be the corresponding derivation of $\FK{\GD}$ in $1$. It restricts to a derivation of 
$\FK{\GD}_{int}\subseteq \FK{\GD}$ in $1$. By Theorem \ref{Deraosrf} there exist elements $n_-\in\n_f^-$, $h\in\h$, and $n_+\in\n_f^+$  such that
\begin{eqnarray}\label{comp}
    \phi(x_V v)=\ti{\delta}_x(f_{\phi v}\res{\widehat{G}})=\phi(n_- v)+\phi(h v)+\phi(n_+ v)
\end{eqnarray}
for all objects $V$ of ${\mathcal O}_{int}$, $v\in V$, and $\phi\in V^{(*)}$.

Let $\La\in P^{++}$ and $v_\La\in L(\La)_\La\setminus\{0\}$. By part (b) we find that there exists elements $\ti{n}_-\in \n^-$, $\ti{h}\in\h$, such that
\begin{eqnarray*}
         \phi(x_{L(\La)} v_\La)=  \phi(\ti{n}_- v_\La)+\phi(\ti{h} v_\La)
\end{eqnarray*}
for all $\phi\in L(\La)^{(*)}$. With (\ref{comp}) we get
\begin{eqnarray*}
   \phi(n_- v_\La)+\La(h) \phi(v_\La)=\phi(\ti{n}_- v_\La)+\La(\ti{h})\phi( v_\La)\\
\end{eqnarray*} 
for all $\phi\in L(\La)^{(*)}$. By choosing $\phi\in (L(\La)^{(*)})_\La\setminus\{0\}$ we find
$\La(\ti{h})=\La(h)$. It follows
\begin{eqnarray*}
   \phi(\ti{n}_- v_\La)=\phi(n_- v_\La)
\end{eqnarray*} 
for all $\phi\in L(\La)^{(*)}$. Therefore $\ti{n}_- v_\La = n_- v_\La$. By (c)
we find $n_-=\ti{n}_-\in\n^-$. Inserting in (\ref{comp}) it follows $x=n_- + h + n_+\in\g_f$.
\qed

%%%%%%%%%%%%%%%%%%%%%%%%%%%%%%%%%%%%%%%%%%%%%%%%%%%%%%%%%%%%%%%%%%%%%%%%%%%%%%%%%%%%%%%%%%%%%%%%%%%%%%%%%%%%%%%%%%%%%%%%%%%%%%%%%%%%%%%%%%%%%%%
%
%
%

%
\end{document}